\documentclass[runningheads]{llncs}
\usepackage[T1]{fontenc}
\usepackage{graphicx}
\usepackage{hyperref}       
\usepackage{url}            

\usepackage{amsmath, amssymb, mathtools, amsfonts}
\usepackage{subcaption, array, booktabs}
\usepackage[capitalize,noabbrev]{cleveref}
\usepackage{color}

\usepackage{algorithm, algorithmic}

\newcommand\set[1]{\left\{ #1 \right\}}

\DeclareMathOperator*{\argmax}{arg\,max}

\DeclareMathOperator{\E}{\mathbb{E}}

\begin{document}
\title{Reinforcement Learning for Molecular Dynamics Optimization: A Stochastic Pontryagin Maximum Principle Approach}
%
%
\author{
Chandrajit Bajaj \inst{1,2}
\and Minh Nguyen \inst{1,2}
\and Conrad Li \inst{1}
}
%
%
\institute{
University of Texas at Austin
\and Equal contribution
}
\maketitle              

\begin{abstract}
In this paper, we present a novel reinforcement learning framework designed to optimize molecular dynamics by focusing on the entire trajectory rather than just the final molecular configuration. Leveraging a stochastic version of Pontryagin’s Maximum Principle (PMP) and Soft Actor-Critic (SAC) algorithm, our framework effectively explores non-convex molecular energy landscapes, escaping local minima to stabilize in low-energy states. Our approach operates in continuous state and action spaces without relying on labeled data, making it applicable to a wide range of molecular systems. Through extensive experimentation on six distinct molecules, including Bradykinin and Oxytocin, we demonstrate competitive performance against other unsupervised physics-based methods, such as the Greedy and NEMO-based algorithms. Our method's adaptability and focus on dynamic trajectory optimization make it suitable for applications in areas such as drug discovery and molecular design.
\end{abstract}
\section{Introduction}
\noindent Molecular dynamics involve the study of how molecules move and interact over time, and it plays a crucial role in fields such as drug discovery, materials science, and biophysics. A molecule's properties are determined not just by its final stable configuration but by the entire trajectory it follows during its evolution. Optimizing these trajectories to achieve desired molecular behaviors, such as minimizing potential energy or maximizing binding affinity, presents significant computational challenges. The inherent non-convexity of molecular energy landscapes, compounded by numerous local minima, poses significant challenges in identifying globally optimal solutions.\\

\noindent Historically, many computational methods, ranging from combinatorial algorithms to deep learning models, have been developed to solve specific molecular dynamics problems. These methods typically focus on finding optimal end-state configurations, such as protein folding or side-chain packing, and often rely on predefined datasets or libraries (e.g., rotamer libraries for protein structures). However, they struggle to explore dynamic molecular trajectories and frequently get trapped in local minima.\\

\noindent \textbf{Classical approaches}: Many classical methods \cite{CKS04,D02,DC97,DK93,DK94,LWRR00,MDL97,MIS87,PR87,TEHL91,X05,XB06} for molecular dynamics  treat optimization as a combinatorial search problem. Given a molecular structure, these methods aim to determine optimal configurations that minimize system energy. However, they are typically constrained to discrete state or action spaces and rely on heuristics or fixed libraries of configurations, making them prone to local minima and less effective for exploring large, continuous search spaces. Alternative techniques such as Markov Chain Monte Carlo\cite{HS91}, cyclical search \cite{DK93,XH01}, spatial restraint 
satisfaction \cite{SB93}, semi-definite programming \cite{CKS04,KCS05,XB06} also struggle to escape local minima.\\

\noindent \textbf{Deep learning approaches}: Recent advancements in deep learning have introduced models \cite{CNN4SCP,SIDEpro,Ingraham:2019,RGN} capable of predicting molecular properties, such as protein folding and binding affinities. These models, including prominent ones like AlphaFold \cite{AF,AF2}, leverage large labeled datasets and focus on final molecular configurations. While these approaches achieve impressive results in predicting final configurations, they often overlook the full dynamic trajectory and are limited by the availability of labeled data \cite{ProteinNet}, restricting their applicability to more general molecular systems.\\

\noindent \textbf{Reinforcement learning approach}: In contrast to classical and deep learning methods, our reinforcement learning (RL) framework addresses the broader challenge of optimizing general molecular dynamics, focusing on the entire trajectory rather than just the endpoint configuration. Our approach captures the entire molecular trajectory, enabling fine-grained optimization of the transitional path between states. Moreover, our framework actively escapes local minima, ensuring more effective exploration of complex energy landscapes. Additionally, our RL approach operates in an unsupervised setting, meaning it does not require labeled training datasets, unlike deep learning models. This flexibility makes our method broadly applicable to various molecular systems, without relying on predefined libraries or supervised data.\\

\noindent Our method uses a reinforcement learning framework that integrates Soft Actor-Critic (SAC) for policy optimization and the Pontryagin Maximum Principle (PMP) to capture molecular dynamics, allowing for continuous control of molecular trajectories. This approach enables us to model molecular dynamics in a way that captures the full continuous energy and force-field interactions. Furthermore, the RL framework efficiently explores the search space, escaping local minima through strategic exploration and exploitation of the molecular energy landscape. Our key contributions include:
\begin{itemize}
\item A reinforcement learning framework based on a stochastic parameterized Hamiltonian formulation of the Pontryagin Maximum Principle (PMP), which optimizes not only the final molecular configuration but also the entire dynamic trajectory. This provides a more comprehensive solution for molecular dynamics, enabling the model to account for continuous energy flows and interactions over time.
\item Unsupervised learning capabilities that leverage physics-based principles, such as stochastic Hamiltonian dynamics, eliminating the need for labeled datasets. This allows the method to generalize effectively across diverse molecular systems without relying on predefined data, making it adaptable to various molecular environments.
\item Efficient exploration strategies embedded within the reinforcement learning framework, which incorporate optimal control theory through stochastic PMP to guide the system out of local minima. This ensures robust performance in navigating complex, non-convex energy landscapes and discovering globally optimal molecular configurations.
\end{itemize}

The rest of the paper is structured as follows: \cref{sec:intro_scp} provides a detailed explanation of the molecular energy-minimization formulation and outlines the application of the soft actor-critic reinforcement learning algorithm to this problem. In \cref{sec:stochastic_pmp}, we introduce the stochastic version of Pontryagin’s Maximum Principle (PMP) and demonstrate its application in selecting an appropriate reward function for the reinforcement learning algorithm. This reward function effectively balances exploration, exploitation, and stabilization during training. Finally, in \cref{sec:experiment_scp}, we present the competitive performance of our approach in learning the molecular dynamics of several essential molecules. Along with fundamental physics-based benchmarks, we compare our method to the unsupervised counterpart of the state-of-the-art physics-based method NEMO \cite{Ingraham:2019}, which we refer to as NEMO-based.
\section{Problem statement and reinforcement learning formulation}\label{sec:intro_scp}
\subsection{Problem statement}
In the context of molecular energy-minimization, the primary learning problem is to predict the optimal molecular configuration that minimizes the system’s potential energy. To achieve this, we account for key molecular configuration variables, including the 3D coordinate $x_i$, bond lengths $b_i$, bond angles $a_i$, and its dihedral angles $d_i$ of the $i$th atom with $i \in \overline{1, N}$. Here $N$ is the number of atoms in the molecule. The following equation establishes a relationship between the 3D coordinates $x = (x_i)_{i=1}^N$ and the dihedral angles $d = (d_i)_{i=1}^N$: 
\begin{equation}\label{eq:inner_transform}
x_i = x_{i-1}  + b_i \begin{bmatrix}u_{i-1} & n_{i-1} \times u_{i-1}  & n_{i-1}\\ \end{bmatrix} \begin{bmatrix} \cos(\pi-a_i) \\ \sin(\pi-a_i)\cos(d_i) \\ \sin(\pi-a_i)\sin d_i \end{bmatrix}
\end{equation}
where $u_i$ is the unit vector from $x_{i-1}$ to $x_i$ and $n_i$ is the unit vector normal to each bond plane.\\

\noindent While the energy function $U$ is often more naturally expressed in terms of the 3D coordinates, the dynamics are better represented in terms of the dihedral angles. To simplify the optimization, we choose to work with the dihedral angles $d$. The primary task then is to find the optimal set of dihedral angles $d$ that minimizes the system's total potential energy. This can be formulated as minimizing $U = U(x) = U(H(d)) = \tilde{U}(d)$, where $U$ is the potential energy function depending on 3D coordinate $x$, and $H$ is the transformation defined by \cref{eq:inner_transform} between 3D coordinates and dihedral angles.\\

\noindent In particular, the potential energy $U$ typically consists of two components, bonded energy and non-bonded energy, that are expressed as:
\begin{align}
U &= U_{total} = U_{bonded} + U_{non-bonded} \nonumber\\
&=\bigg(\sum_{bonds} K_r(r-r_0)^2 + \sum_{angles} K_{\theta}(\theta-\theta_0) + \sum_{rotamers} f(rotamer) \bigg) \\
&+ \sum_{dihedrals} f(dihedral) + \left(\sum_{i < j} \frac{q_iq_j}{\epsilon r_{ij}} + \sum_{i < j} \left(\frac{A_{ij}}{r_{ij}^{12}} - \frac{B_{ij}}{r_{ij}^6}\right) \right) 
\end{align}
where the first group of terms describes bonded interactions (stretching, bending, and torsion), and the second group captures non-bonded interactions (electrostatics and van der Waals forces, modeled by the Lennard-Jones potential).\\

\noindent In our approach, rather than solely focusing on finding the optimal configuration $c^*$, we aim to optimize the entire trajectory of the molecule as it transitions from an initial configuration $c^0$ to the final configuration $c^*$. This trajectory-based optimization allows us to model and capture the intermediate states and dynamics of the system, ensuring that the molecule avoids local minima and follows a path that effectively minimizes the potential energy over time. The reinforcement learning framework is used to learn and optimize this trajectory, providing a more complete solution to the molecular energy-minimization problem.

\subsection{Reinforcement learning formulation}
\noindent Building on the previous section, molecular dynamics optimization can be framed as a reinforcement learning problem, where the state is defined by the dihedral angles $d$. At each time step $t$, the molecular configuration evolves based on the change in dihedral angles, described by their velocities $\dot{d}_t$ at time $t$. Thus, the state $s_t$ is the current set of dihedral angles $d_t$, and the action $a_t$ is the adjustment of these angles at time $t$.\\

\noindent The system transitions from state $s_t$ to $s_{t+1}$ according to the discretized transition dynamics $s_{t+1} = s_t + \epsilon a_t$. Here we select the constant step size $\epsilon = 0.3$ for all experiments shown in \cref{sec:experiment_scp}. The trajectory ends when either a state violates physically plausible ranges or a maximum number of steps has been reached. The reward function $r_t = r(s_t, a_t)$ depends on the potential energy $U$ and how it evolves with the action $a_t$. This reward is designed to guide the system toward minimizing energy over the trajectory. An explicit formulation for the reward function is provided in Section 3. 

\subsection{Soft actor-critic}
In this paper, we will use the soft actor-critic as our main reinforcement learning algorithm. The goal of reinforcement learning in this context is to learn the optimal stochastic policy $\pi(a_t|s_t)$ which maximizes the expected sum of rewards over time. To enhance exploration and avoid premature convergence to local minima, the soft actor-critic (SAC) method augments this objective by adding an entropy term $\mathcal{H}$:
\begin{equation}
\pi^* = \argmax_{\pi} \sum_t \E_{(s_t, a_t) \sim \rho_{\pi}}[r(s_t, a_t) + \alpha \mathcal{H}(\pi(.|s_t))]
\end{equation}
Next, the model is trained to minimize the soft Bellman residual:
\begin{equation}
J_Q(\phi) = \E_{(s_t, a_t) \sim \mathcal{D}}\bigg[\frac{1}{2}\bigg(Q_{\phi}(s_t, a_t) - (r(s_t, a_t) + \E_{s_{t+1} \sim \rho}[V_{\bar{\phi}}(s_{t+1})])\bigg)^2 \bigg]
\end{equation}

\noindent If the policy is parameterized by a neural network $f_{\theta}$, where actions are generated as $a_t = f_{\theta}(s_t; \epsilon_t)$ (with $\epsilon_t$ being noise sampled from a fixed distribution), the objective function for policy optimization is:
\begin{equation}
J_{\pi}(\theta) = \E_{s_t \sim \mathcal{D}, \epsilon_t \sim \mathcal{N}}\Bigg[\alpha \log \pi_{\theta}(f_{\theta}(s_t; \epsilon_t)|s_t) - Q_{\phi}(s_t, f_{\theta}(s_t; \epsilon_t))\Bigg]
\end{equation}

\noindent In the next section, we discuss how to choose an appropriate reward function to optimize molecular trajectories in the context of energy-minimization tasks.
\section{Reward function via stochastic Pontryagin Maximum Principle}\label{sec:stochastic_pmp}
The reward function is a critical component of our reinforcement learning framework, guiding the system through exploration and subsequent stabilization. Based on the stochastic Pontryagin Maximum Principle (PMP), we designed a reward function that drives exploration through action amplification followed by stabilization after a cutoff time $T_0$. We define the following reward function:
\begin{equation}\label{eqn:reward_formula}
r(s, t) = e^{-\lambda t}(a^2/2 - cU(s)) 1_{\set{t \leq T_0}} - e^{\lambda t}(a^2/2 + cU(s)) 1_{\set{t > T_0}}
\end{equation}
For all experiments in \cref{sec:experiment_scp}, we select $\lambda = \log(1.001) > 0$. Details for hyperparameters $c$ is given in the code link in \cref{sec:conclusion}.

\noindent This reward function modulates the balance between action magnitude $a$ and the system's potential energy $U(s)$, where the parameter $\lambda$ controls the transition between the exploration and stabilization phases, with $T_0$ being the stabilization cutoff time. The exploration phase (for $t \leq T_0$) encourages larger actions to thoroughly explore the state space, while the stabilization phase (for $t > T_0$) dampens the actions, driving the system towards a stable, low-energy configuration. The reward function $r(s_t, a_t)$ is divided into two parts:
\begin{enumerate}
    \item Exploration phase ($t \leq T_0$): The first term encourages exploration by amplifying action magnitudes. This mechanism ensures the system effectively explores the state space, optimizing its trajectory based on the gathered information.
    \item Stabilization phase ($t > T_0$): The second term emphasizes stabilization and refinement by suppressing large actions, guiding the system toward a stable, low-energy configuration.
\end{enumerate}

\noindent The behavior of the system in each phase is characterized by the following dynamics detailed in \cref{lem:dynamic_form} below.

\begin{lemma}\label{lem:dynamic_form}
Let $\Gamma(t)$ be a time-dependent function. For a reward function of the form $r(s, a) = \Gamma(t)(\pm a^2/2 - V(s))$, the resulting optimal dynamics follow the trajectory described by:
\begin{align}
ds_t &= a_t dt + \sigma dW_t \\
da_t &= -\frac{\Gamma'(t)}{\Gamma(t)} a_t dt \mp \nabla V(s_t) dt + K dW_t
\end{align}
where $W_t$ is the white noise, and $K$ is a constant related to the covariance $\sigma$.
\end{lemma}

\noindent We apply \cref{lem:dynamic_form} to two distinct stages of the reward function before and after the cutoff time $T_0$:
\begin{align}
\Gamma_1(t) &= e^{-\lambda t},\ r_1(s_t, a_t) = \Gamma_1(t) (a^2/2 - cU(s)) \\
\Gamma_2(t) &= e^{\lambda t},\ r_2(s_t, a_t) = \Gamma_2(t) (-a^2/2 - cU(s))
\end{align}

\noindent Combining the resulting dynamics from these two stages, we obtain the following final dynamics:
\begin{align}
ds_t &= a_t dt + \sigma dW_t \\
da_t &= \lambda (1_{\set{t \leq T_0}} - \lambda 1_{\set{t > T_0}}) a_t dt - c\nabla U(s_t)(1_{\set{t \leq T_0}} - 1_{\set{t > T_0}})dt + K dW_t
\end{align}

\noindent In this framework, before the cutoff time $T_0$, the action increment is driven by $\lambda a_t - c\nabla U(s_t)$, which not only follows the gradient descent direction to reduce the potential energy but also amplifies the action $a_t$. As a result, larger actions $a_t$ lead to even larger increments, accelerating the system towards its goal. After $T_0$, however, if significant action has already been taken (i.e., $a_t$ is large), the system stabilizes, as $da_t \approx -\lambda a_t$, causing the magnitude of $a_t$ to decrease and resulting in a damping effect.

\begin{remark}
By applying \cref{lem:dynamic_form} to $\Gamma(t) = t^3$, and $r(s, a) = t^3(a^2/2-U(s))$, we recover the continuous dynamics of Nesterov's accelerated gradient method \cite{Su2015-pe}:
\begin{equation}
\ddot{s} + \frac{3}{t}\dot{s} + \nabla U(s) = 0
\end{equation}
Since the exponential factor in our current formulation has a stronger impact on the magnitude of the dynamics, we focus on $\Gamma(t) = e^{\pm \lambda t}$ in this work. Exploring Nesterov-type dynamics will be the focus of future research.
\end{remark}

\subsection{Soft Actor-Critic Algorithm Implementation}
The soft-actor critic algorithm is designed to optimize the policy by maximizing the expected cumulative reward, as described earlier. We integrate the reward function described above into the final algorithm (\cref{alg:SAC_molecular_dynamics}).

\begin{algorithm}
\caption{Main algorithm via Soft Actor-Critic}\label{alg:SAC_molecular_dynamics}
    \begin{algorithmic}[1]
    \STATE \textbf{Input}: Q-function parameters $\phi_1, \phi_2$, empty replay buffer $\mathcal{D}$, maximum number of steps $T_{max}$, hyperparameter $\tau$, and number of iterations $I_{max}$
    \STATE Set target parameters to Q-function parameters $\bar{\phi}_i \gets \phi_i$ for $i \in \set{1, 2}$
    \REPEAT
        \STATE Observe state $s_t=(d_j)_{j=1}^N$.
        \STATE Select action $a_t \sim \pi_{\theta}(.|s_t)$ based on the policy $\pi_{\theta}$
        \STATE Execute $a_t$ to obtain new state $s_{t+1} = (d_j^{(t+1)})_{j=1}^N$
        \STATE Calculate coordinates $(x^{(t)}_i)_{i=1}^N$ and use that to calculate potential  energy $U_t$.
        \STATE Compute the reward function $r(s_t, a_t)$ using the formula \cref{eqn:reward_formula}
        \STATE Update replay buffer $\mathcal{D} \gets \mathcal{D} \cup \set{(s_t, a_t, r(s_t, a_t), s_{t+1}}$
        \STATE If $t \geq T_{max}$ or $s_{t+1}$ violate physical range constraint, $s_t$ is considered terminal state. In this case, we reset the state by uniformly randomly sampled dihedral angles.
        \IF {update step}
            \STATE Use randomly sample batch $\mathcal{B}$ from $\mathcal{D}$ to update $Q$-function and policy.
            \STATE Soft update target network $\bar{\phi_i} \gets \tau \phi_i + (1-\tau)\bar{\phi}_i$ for $i \in \set{1, 2}$
        \ENDIF
    \UNTIL{the number of iterations exceeds $I_{max}$}
    \end{algorithmic}
\end{algorithm}

\begin{figure}[!htb]
\minipage{0.22\textwidth}
  \includegraphics[width=\linewidth]{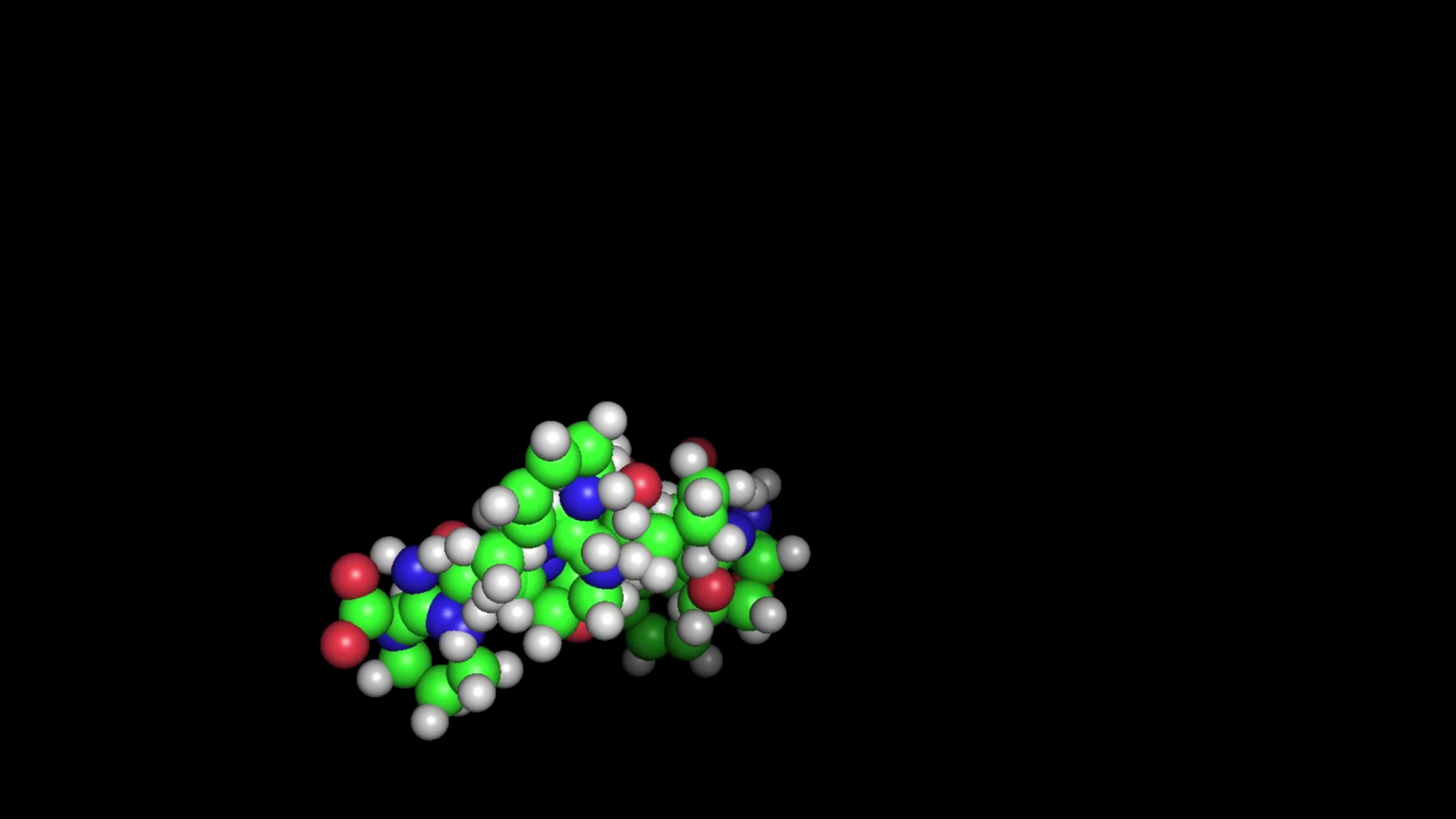}
\endminipage\hfill
\minipage{0.22\textwidth}
  \includegraphics[width=\linewidth]{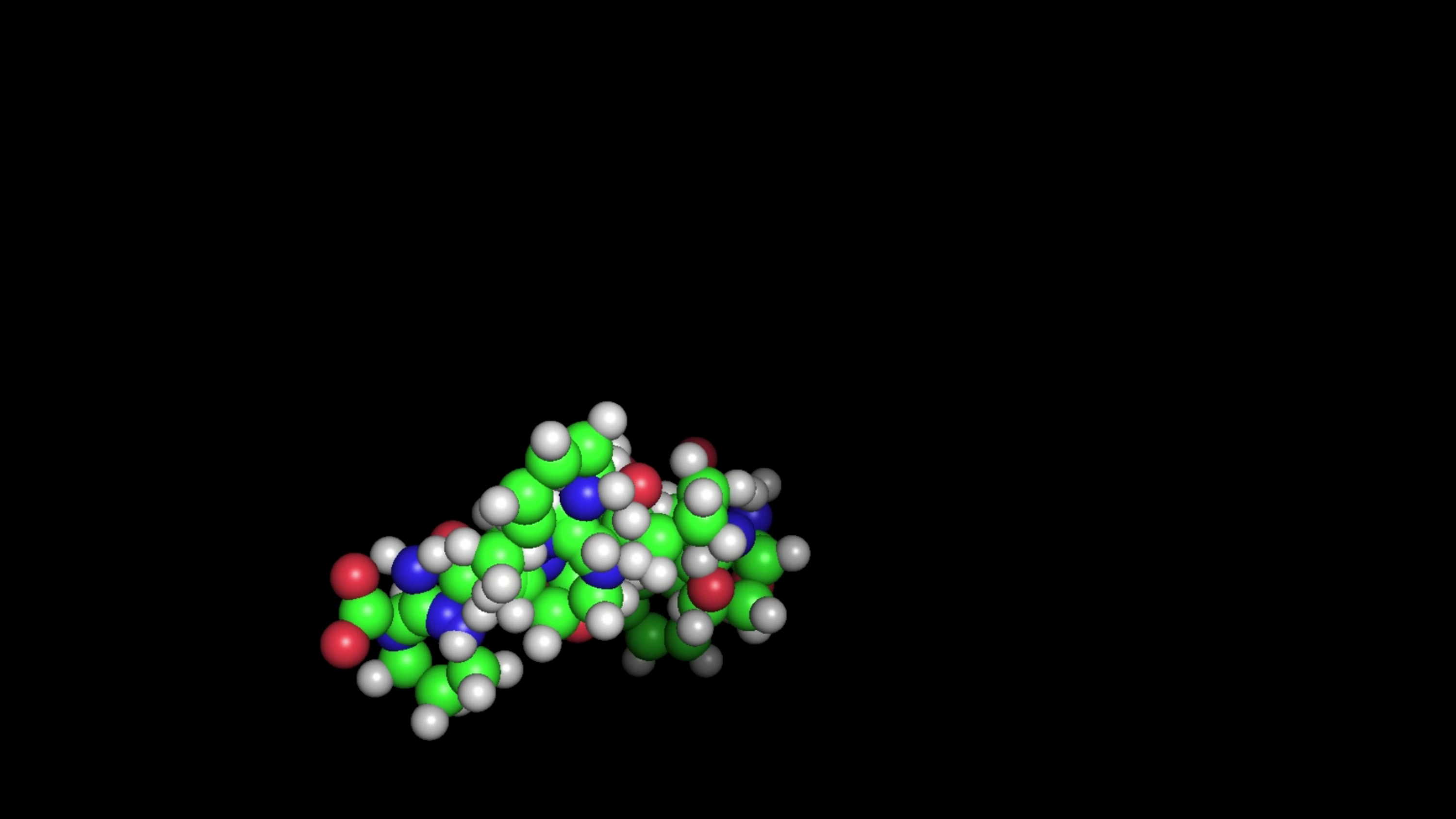}
\endminipage\hfill
\minipage{0.22\textwidth}%
  \includegraphics[width=\linewidth]{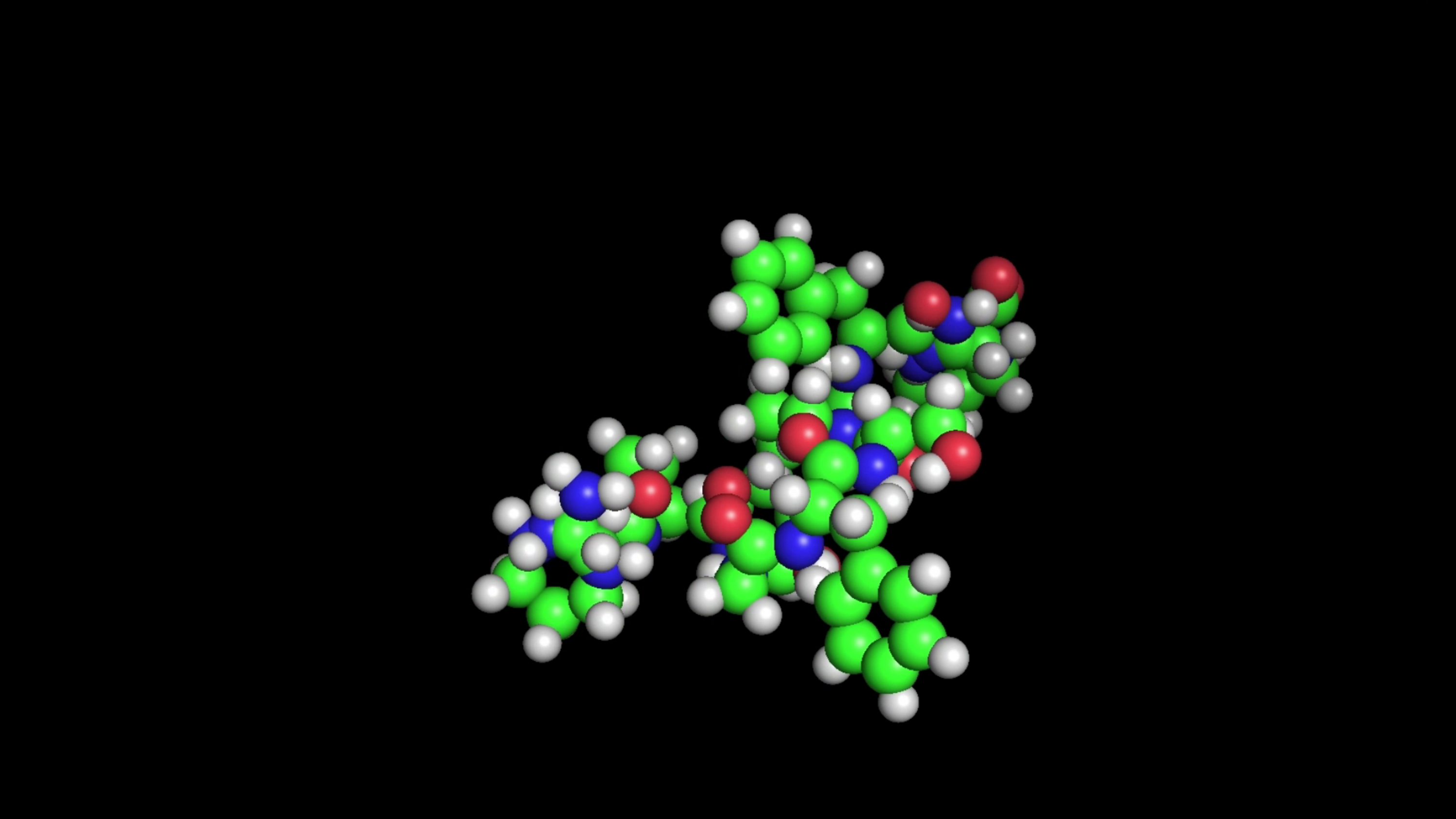}
\endminipage\hfill
\minipage{0.22\textwidth}%
  \includegraphics[width=\linewidth]{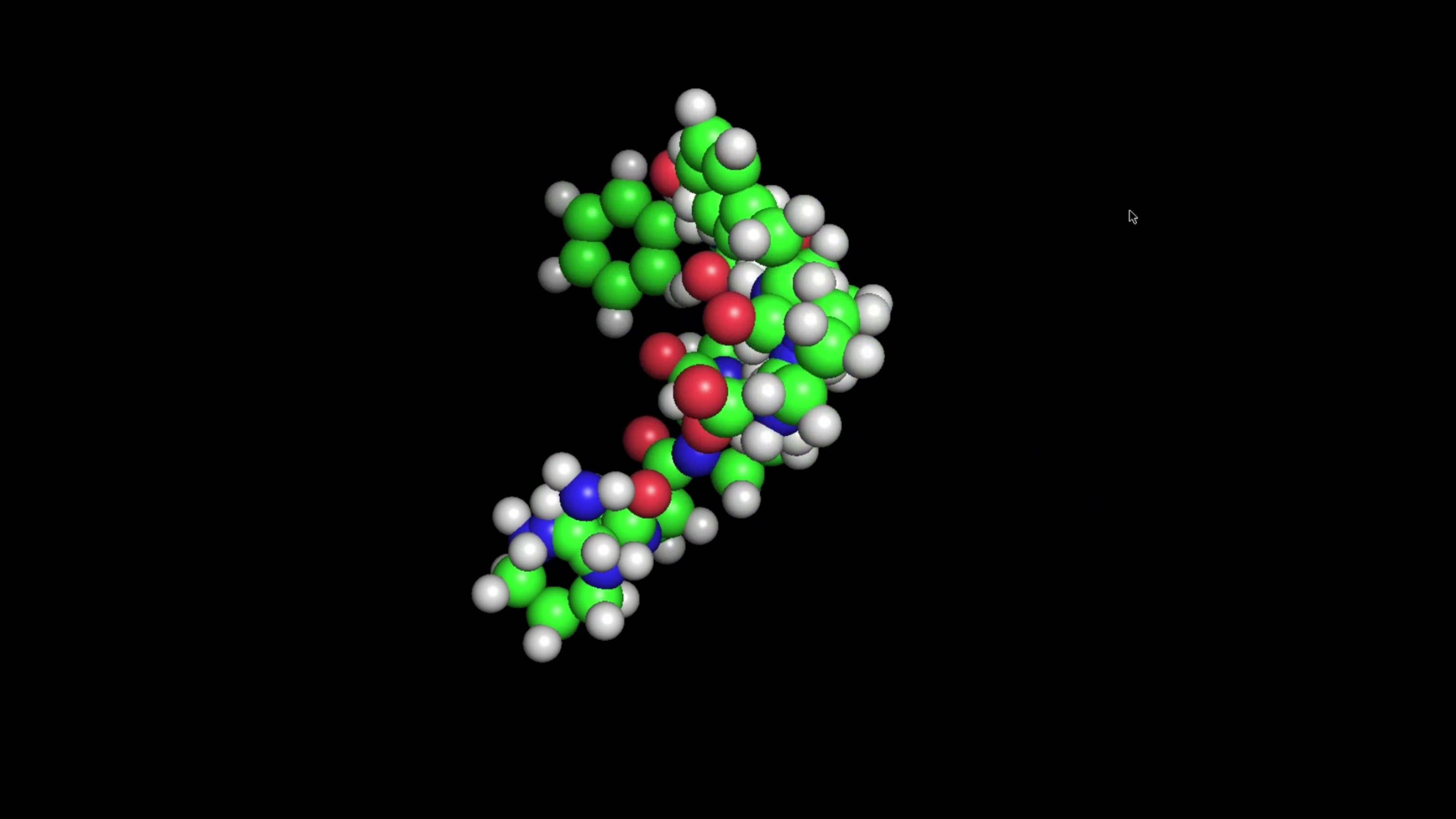}
\endminipage
\caption{Snapshots from a learned episode for Bradykinin molecule}
\label{fig:visualize_bradykinin}
\end{figure}
\section{Experiments}\label{sec:experiment_scp}
In this section, we evaluate our reinforcement learning-based approach for learning optimal molecular dynamics across six different molecules: Bradykinin, CLN025, Met-enkephalin, Oxytocin, Substance-P, and Vasopressin. Our objective is to minimize potential energy while optimizing the full molecular trajectory, rather than focusing solely on the final configuration.

\begin{figure}[!htb]
\minipage{0.22\textwidth}
  \includegraphics[width=\linewidth]{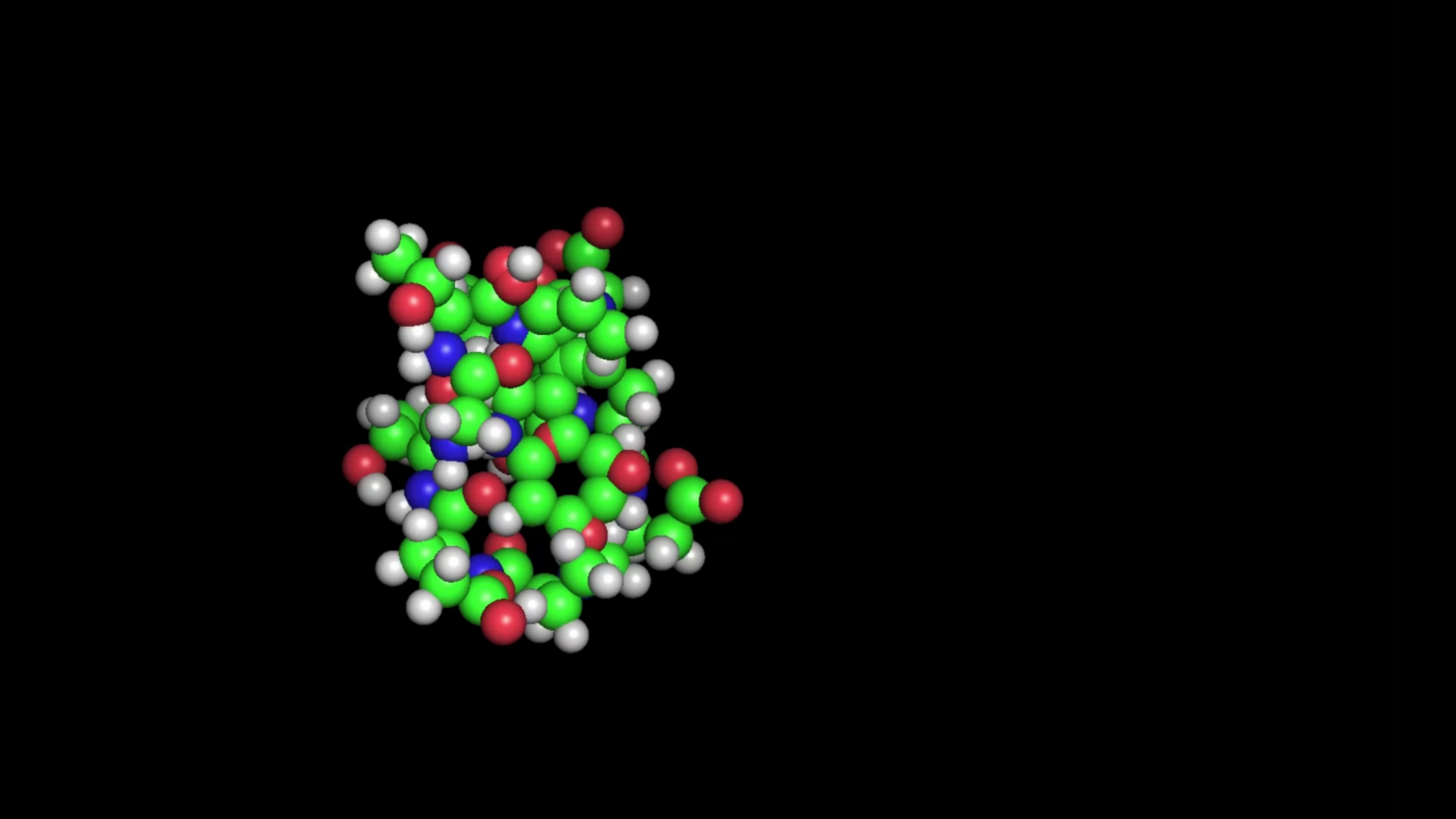}
\endminipage\hfill
\minipage{0.22\textwidth}
  \includegraphics[width=\linewidth]{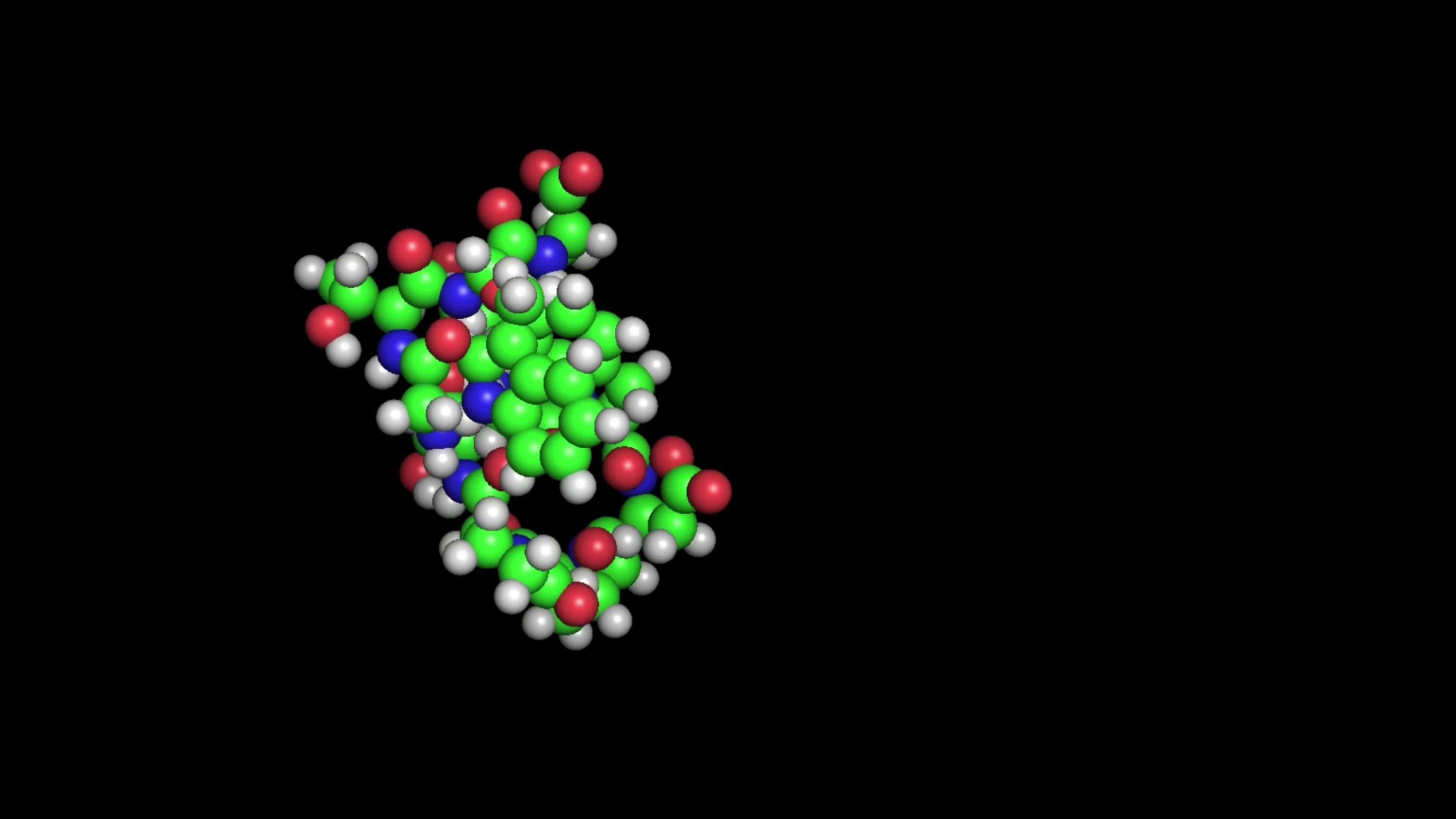}
\endminipage\hfill
\minipage{0.22\textwidth}%
  \includegraphics[width=\linewidth]{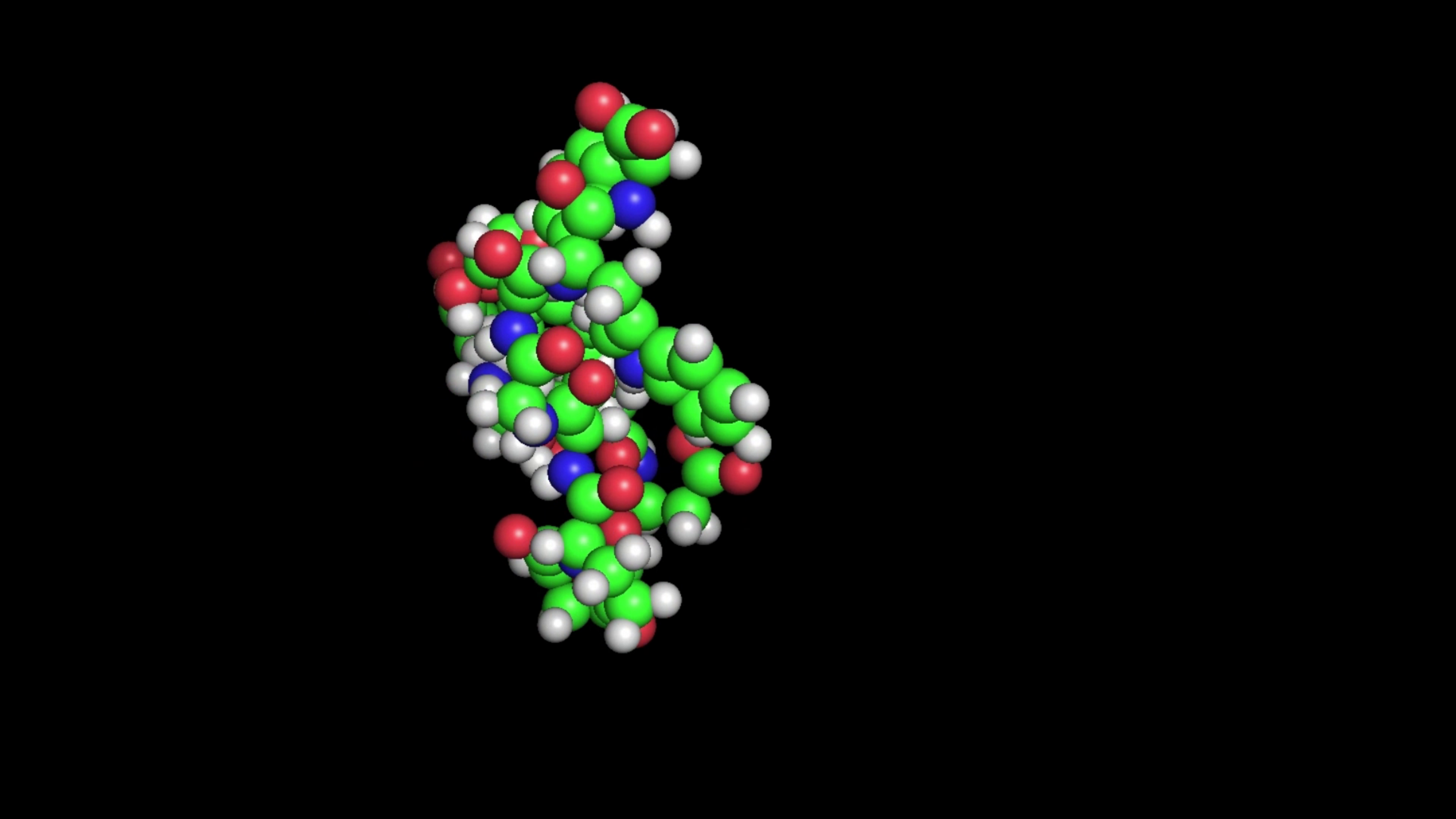}
\endminipage\hfill
\minipage{0.22\textwidth}%
  \includegraphics[width=\linewidth]{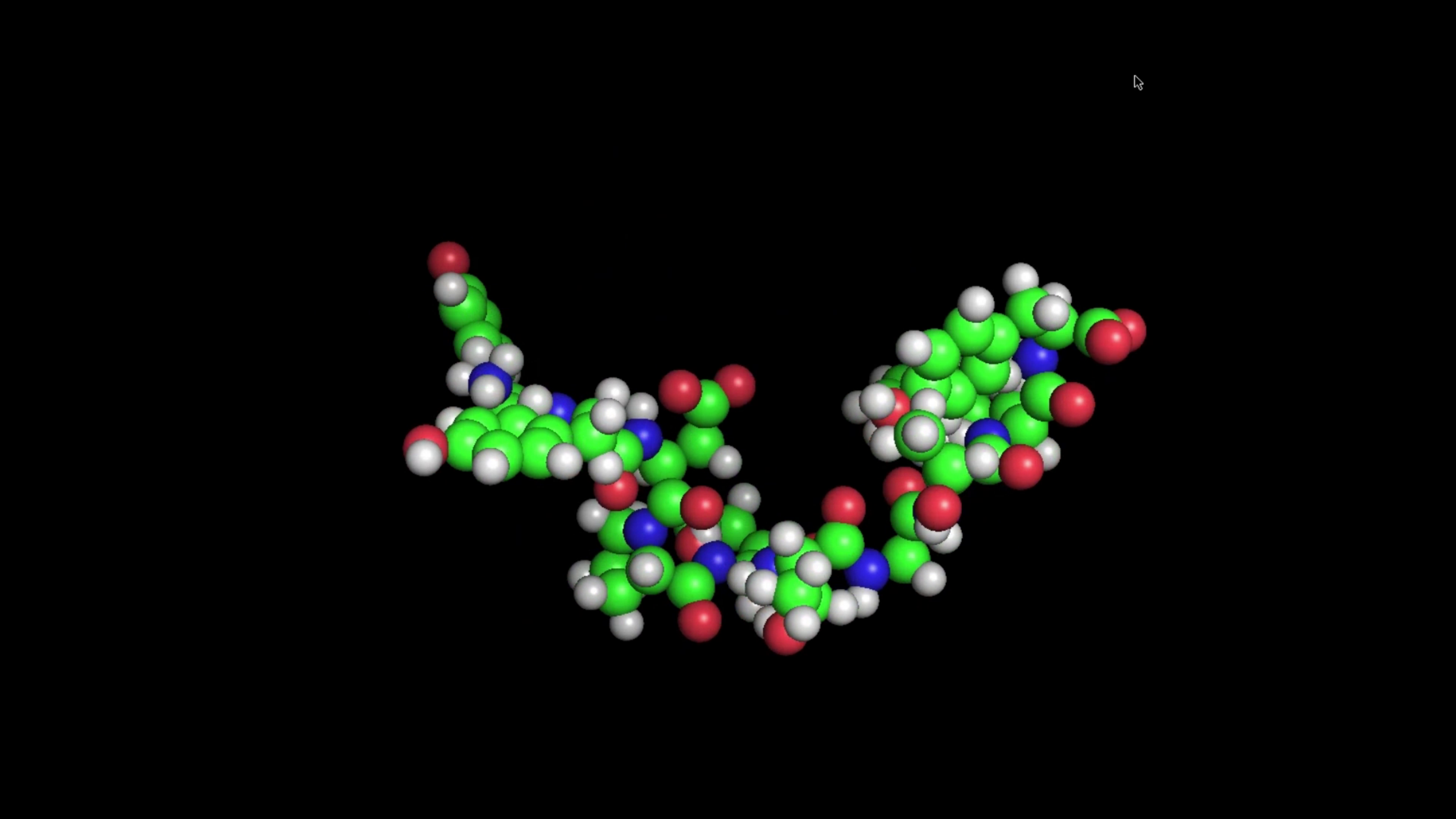}
\endminipage
\caption{Snapshots from a learned episode for CLN025 molecule}\label{fig:visualize_cln025}
\end{figure}

We implement our reinforcement learning-based model as outlined in \cref{alg:SAC_molecular_dynamics}. The policy network takes into account the dihedral angles $d$ along with the fixed configuration feature $c=(X, A)$, where $X$ is the chemical and electrostatic feature vectors, and $A$ is the bond adjacency matrix. These features remain fixed for a given molecule throughout the optimization process. The feature vectors in $X$ include attributes such as the van der Waals radius, atomic charge, acceptor and donor indicators. The energy function $U$ is computed using Rosetta’s energy scoring library \cite{Alford2017}.\\

\begin{figure}[!htb]
\minipage{0.22\textwidth}
  \includegraphics[width=\linewidth]{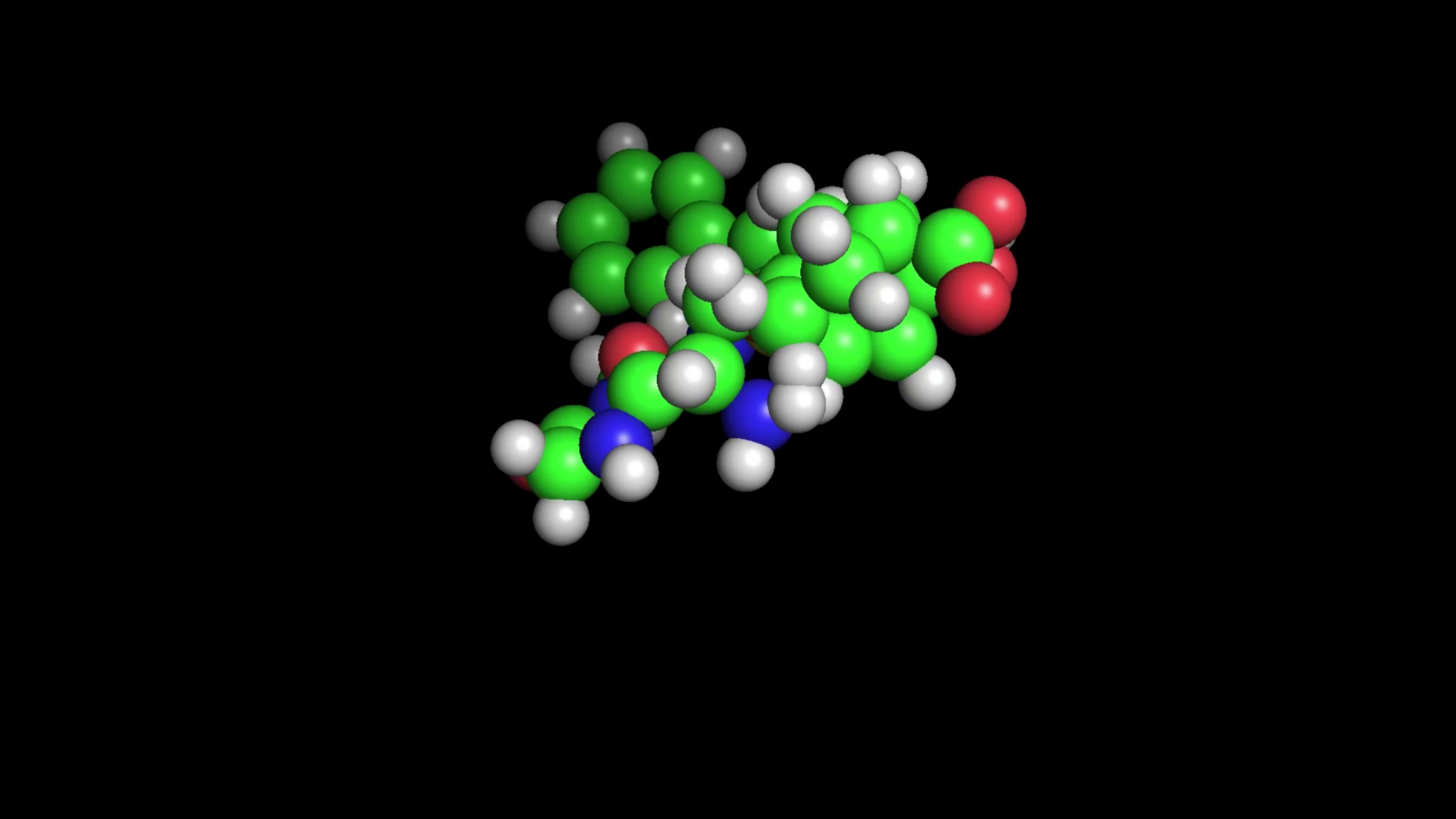}
\endminipage\hfill
\minipage{0.22\textwidth}
  \includegraphics[width=\linewidth]{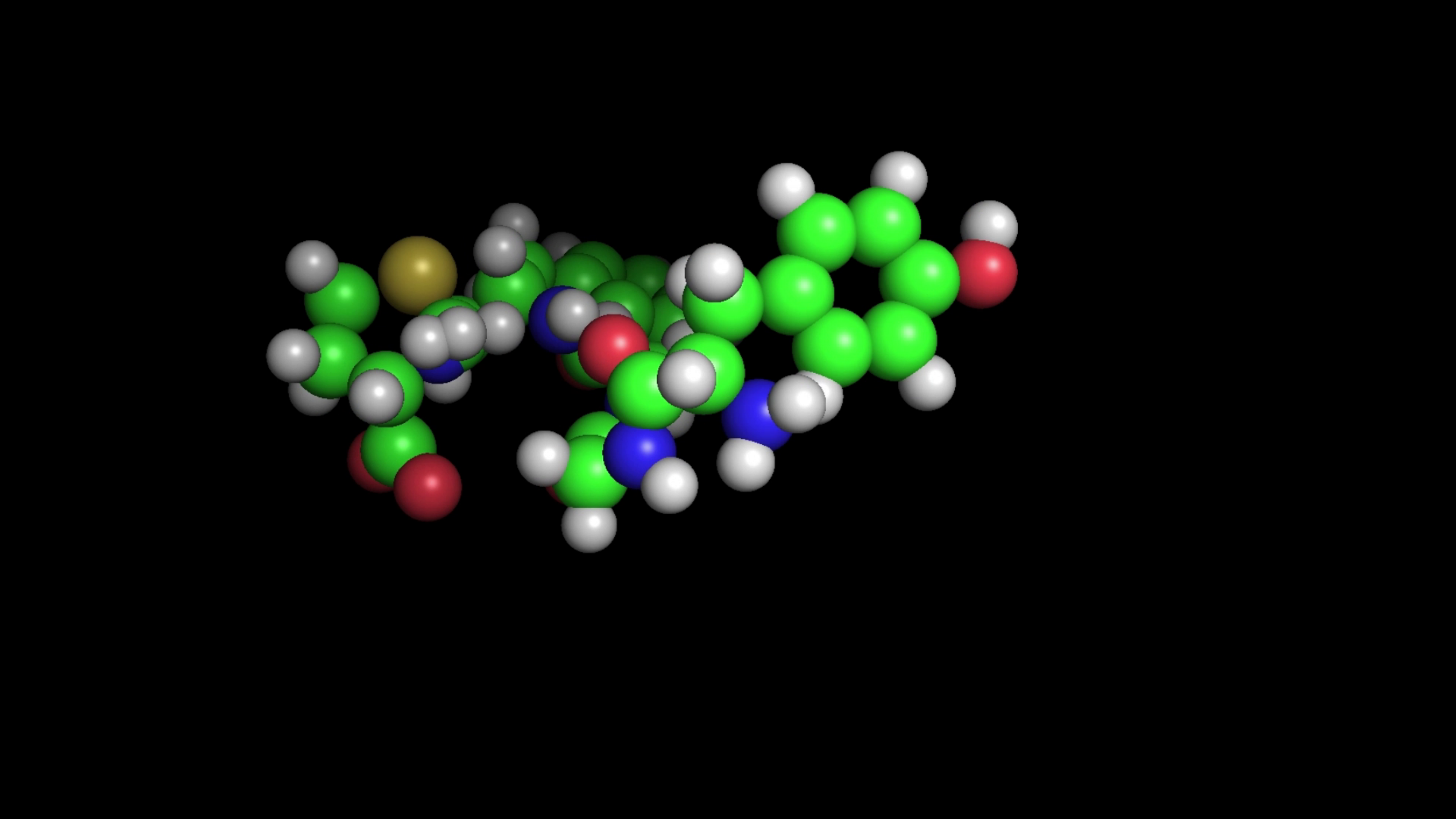}
\endminipage\hfill
\minipage{0.22\textwidth}%
  \includegraphics[width=\linewidth]{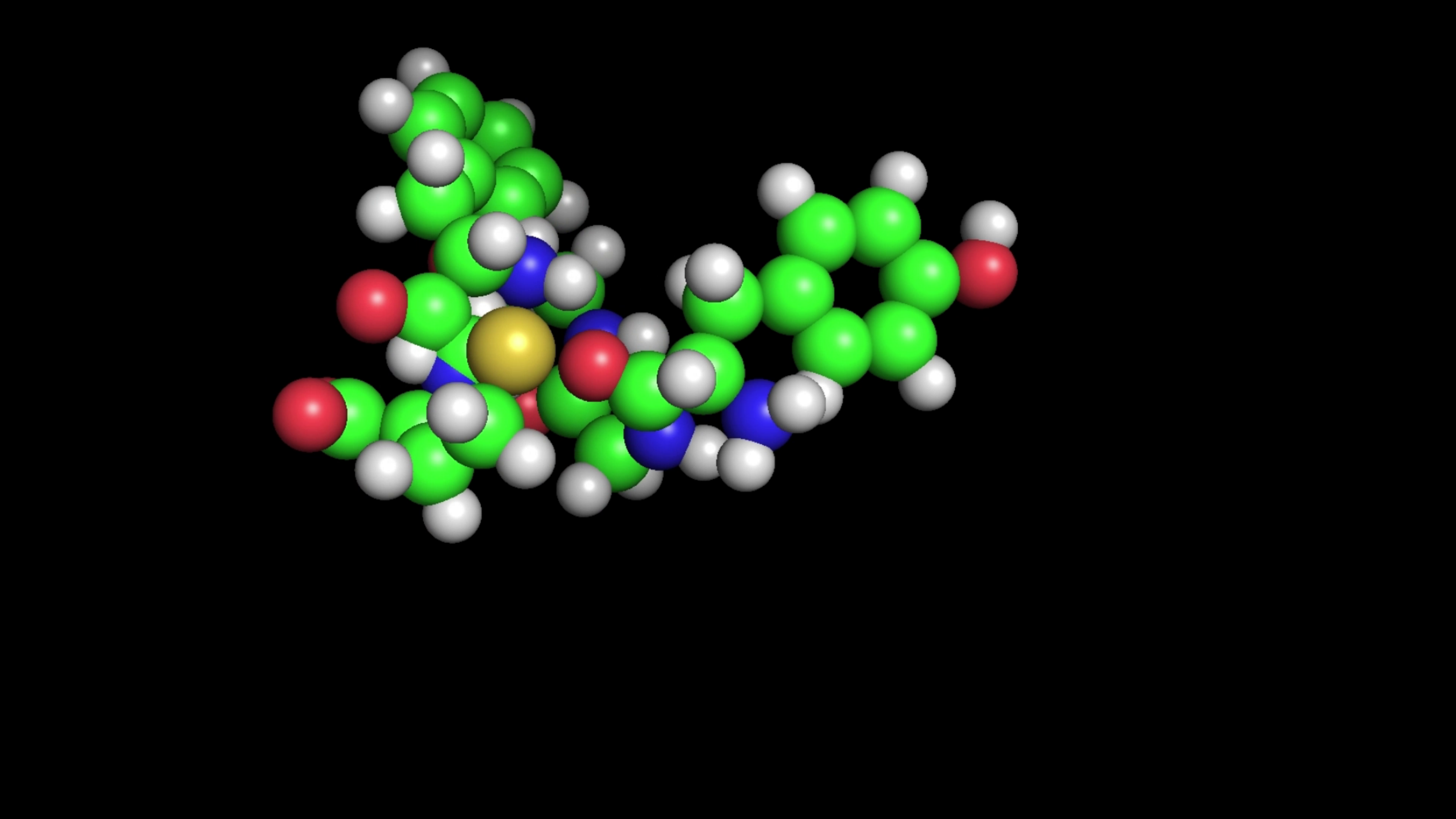}
\endminipage\hfill
\minipage{0.22\textwidth}%
  \includegraphics[width=\linewidth]{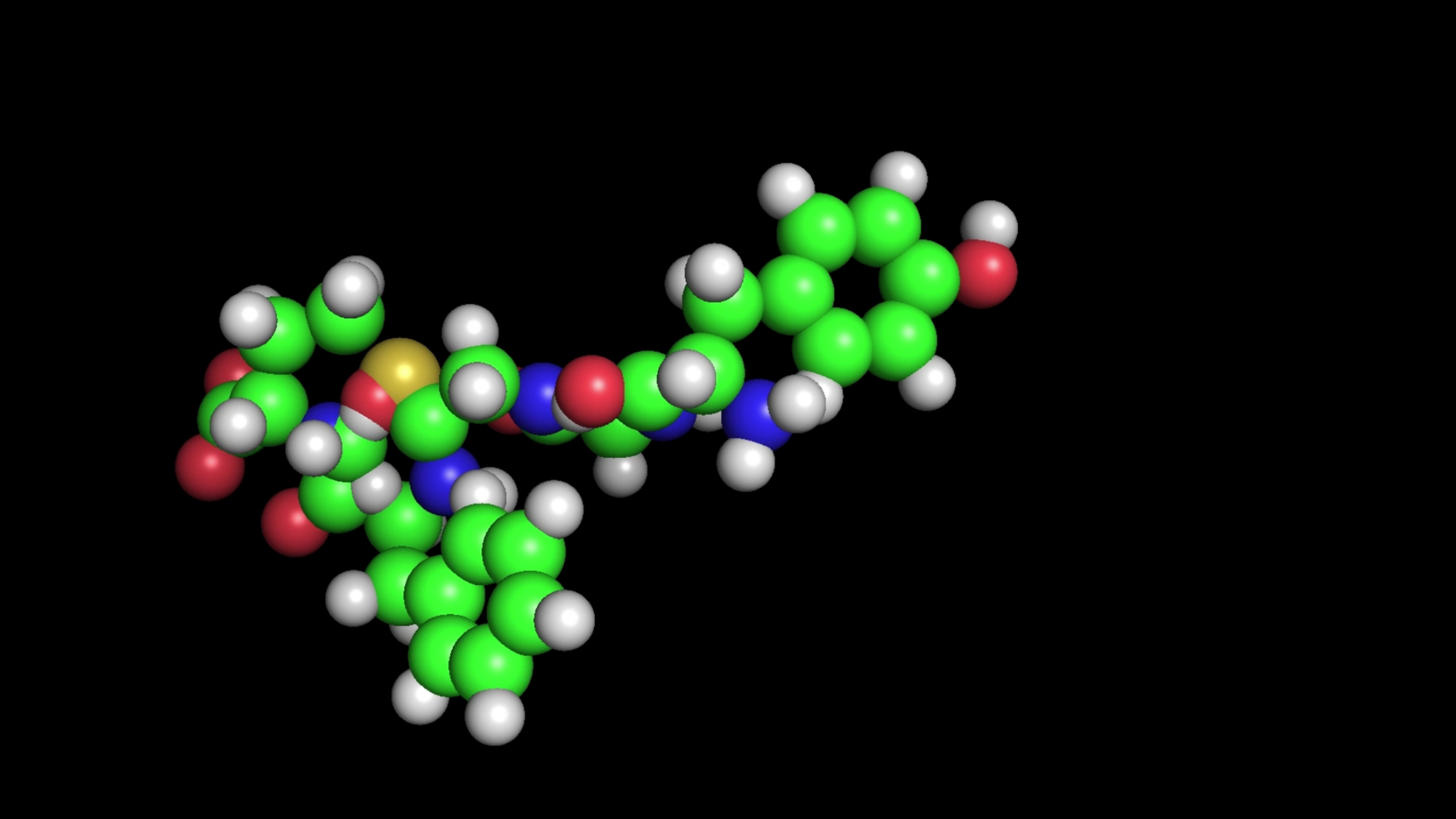}
\endminipage
\caption{Snapshots from a learned episode for Met-enkephalin molecule}\label{fig:visualize_met_enkephalin}
\end{figure}

\begin{table}[htbp]
    \caption{Benchmarking results on 6 different molecules using the \textbf{Ratio} metric}
    \label{table:benchmarks}
    \vspace{3mm}
    \centering
    \begin{tabular}{{c@{\hskip 0.1in}c@{\hskip 0.1in}c@{\hskip 0.1in}c@{\hskip 0.1in}c@{\hskip 0.1in}}}
    \toprule
    & \parbox{60pt}{\centering Greedy} & \parbox{60pt}{\centering Random Hamiltonian} & \parbox{60pt}{\centering NEMO-based} & \parbox{60pt}{\centering \textbf{Ours}} \\
    \midrule
    Bradykinin & 0.15 & 0.08 & 0.75 & {\color{red}0.83} \\
    CLN025 & 0.21 & 0.04 & 0.92 & {\color{red}0.98} \\
    Met-enkephalin & 0.23 & 0.07 & 0.78 & {\color{red}0.84} \\
    Oxytocin & 0.37 & 0.05 & 0.96 & {\color{red}1.0} \\
    Substance-P & 0.19 & 0.11 & 0.66 & {\color{red}0.84} \\
    Vasopressin & 0.08 & 0.02 & {\color{red}0.62} & 0.60 \\
    \bottomrule
    \end{tabular}
\end{table}

\begin{figure}[htb]
  \centering
  \subfloat[Bradykinin]{\includegraphics[scale=0.3]{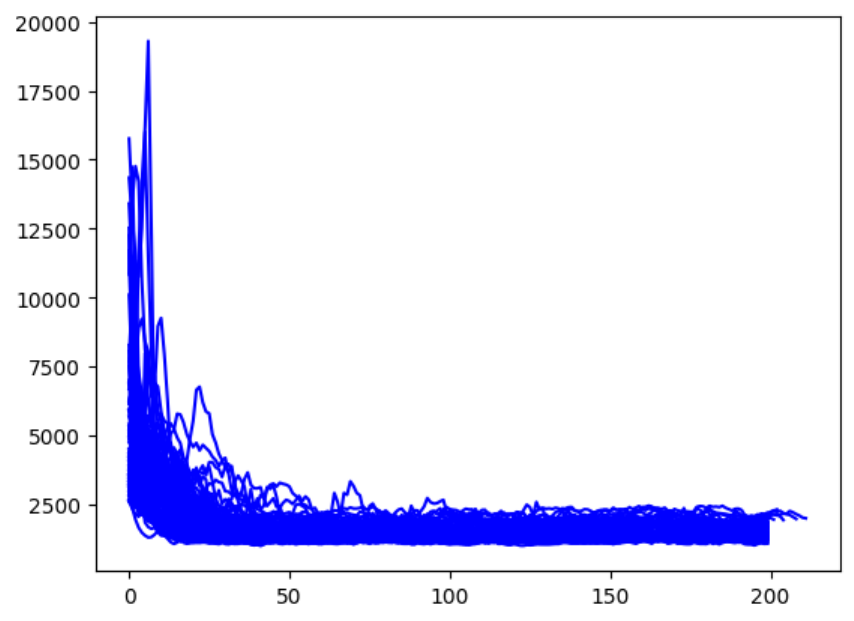}}
  \hfil
  \subfloat[CLN025]{\includegraphics[scale=0.3]{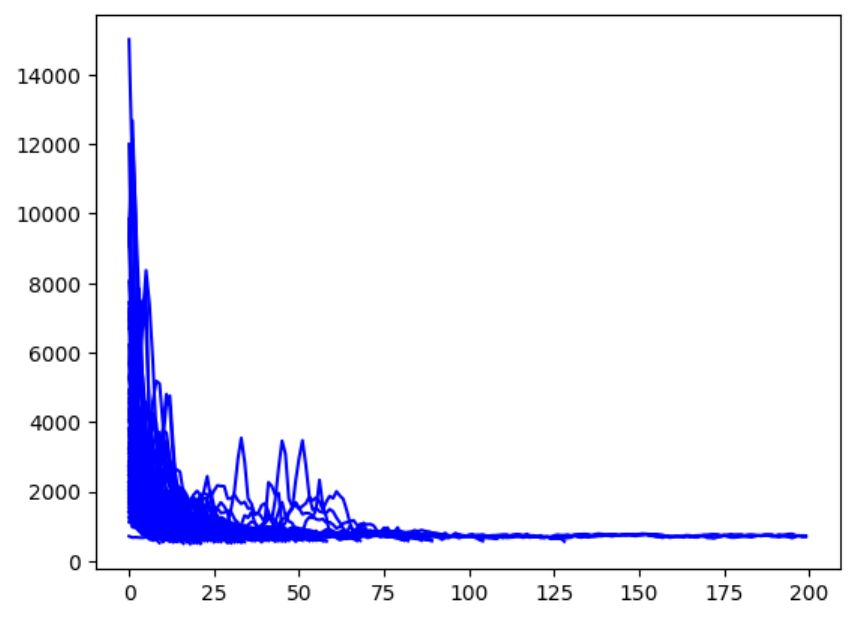}}
  \hfil
  \subfloat[Met-enkephalin]{\includegraphics[scale=0.3]{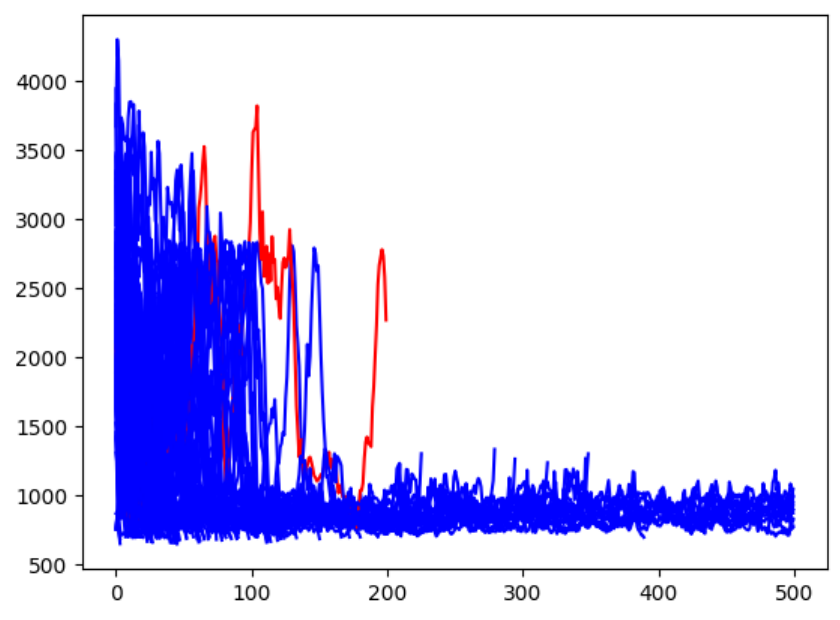}}
  \hfil
  \subfloat[Oxytocin]{\includegraphics[scale=0.3]{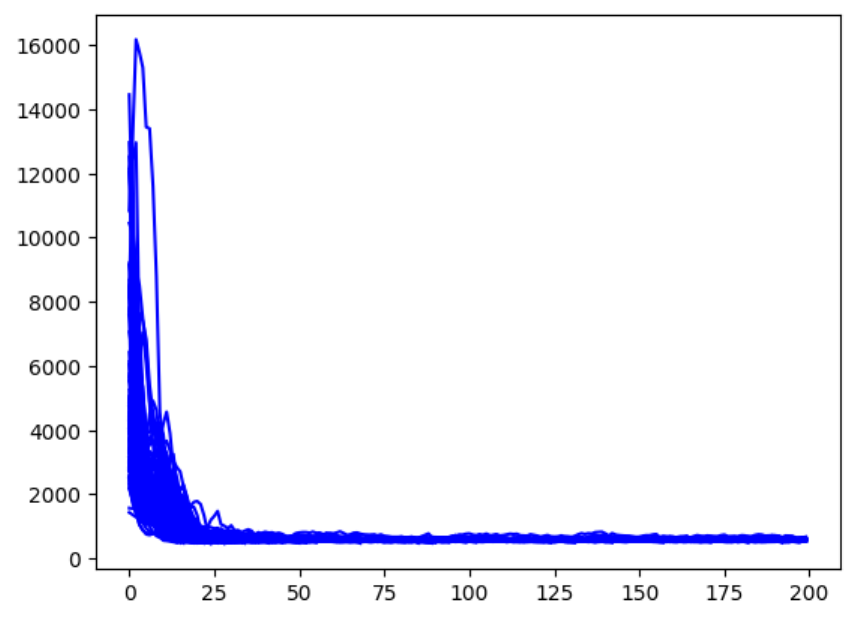}}
  \hfil
  \subfloat[Substance-P]{\includegraphics[scale=0.3]{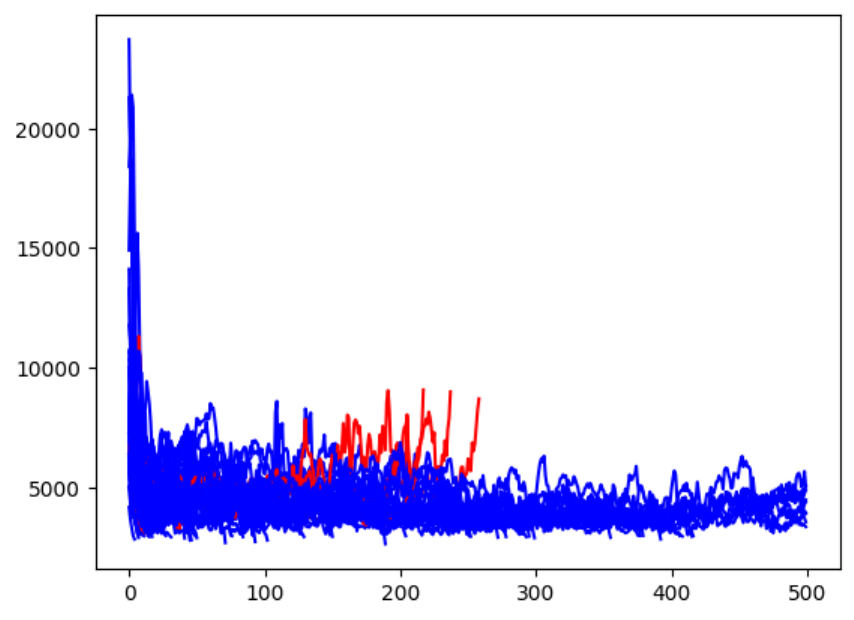}}
  \hfil
  \subfloat[Vasopressin]{\includegraphics[scale=0.3]{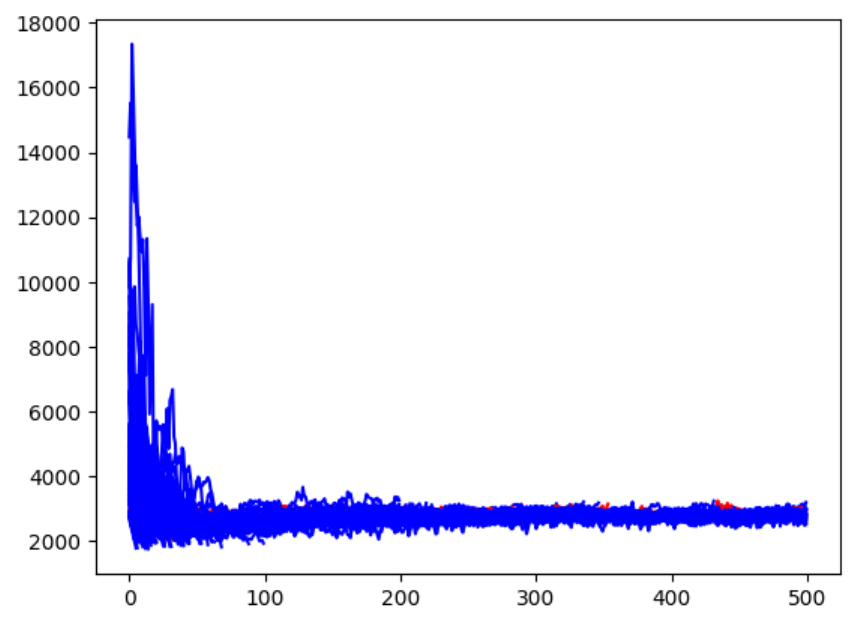}}
  \caption{Energy evolution over time across 100 learned episodes, illustrating the system's stabilization towards low-energy configurations. The $x$-axis represents time steps, while the $y$-axis shows the potential energy values. Red curves are bad learned episodes with high end-energy and large variance.} \label{fig:energy_plot_over_time}
\end{figure}
We compared our approach's performance against three baseline methods:
\begin{enumerate}
    \item Greedy Algorithm: Takes 10 random actions and selects the most energy-reducing action at each time step.
    \item Random Hamiltonian: Uses random updates to the molecular dynamics based on Hamiltonian dual dynamics.
    \item NEMO-based: A neural energy-based model that employs Langevin dynamics to optimize molecular trajectories. This method is an unsupervised version of \cite{Ingraham:2019} and and shares similarities with other state-of-the-art models like AlphaFold \cite{AF}.
\end{enumerate}

\begin{figure}[!htb]
\minipage{0.22\textwidth}
  \includegraphics[width=\linewidth]{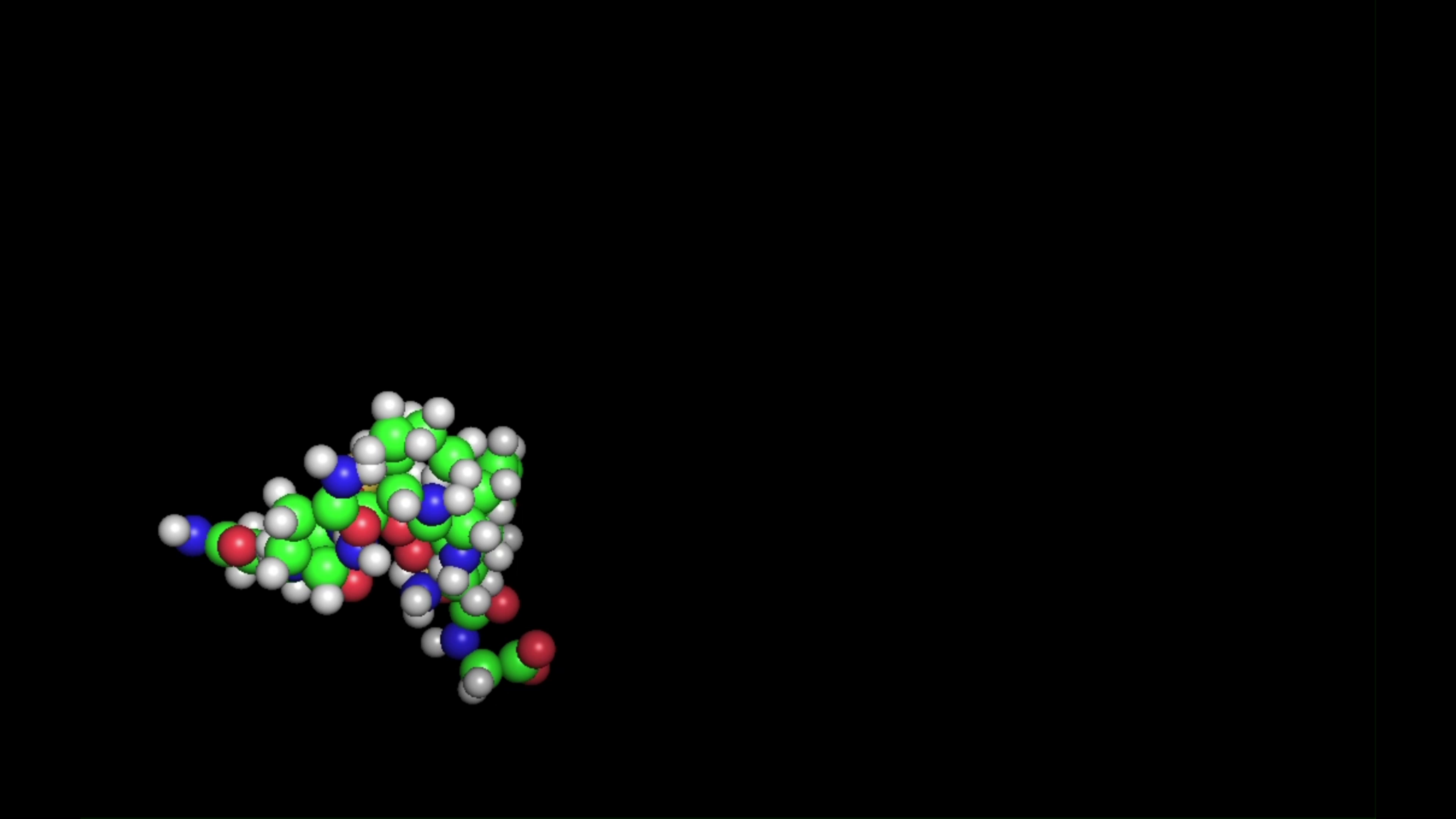}
\endminipage\hfill
\minipage{0.22\textwidth}
  \includegraphics[width=\linewidth]{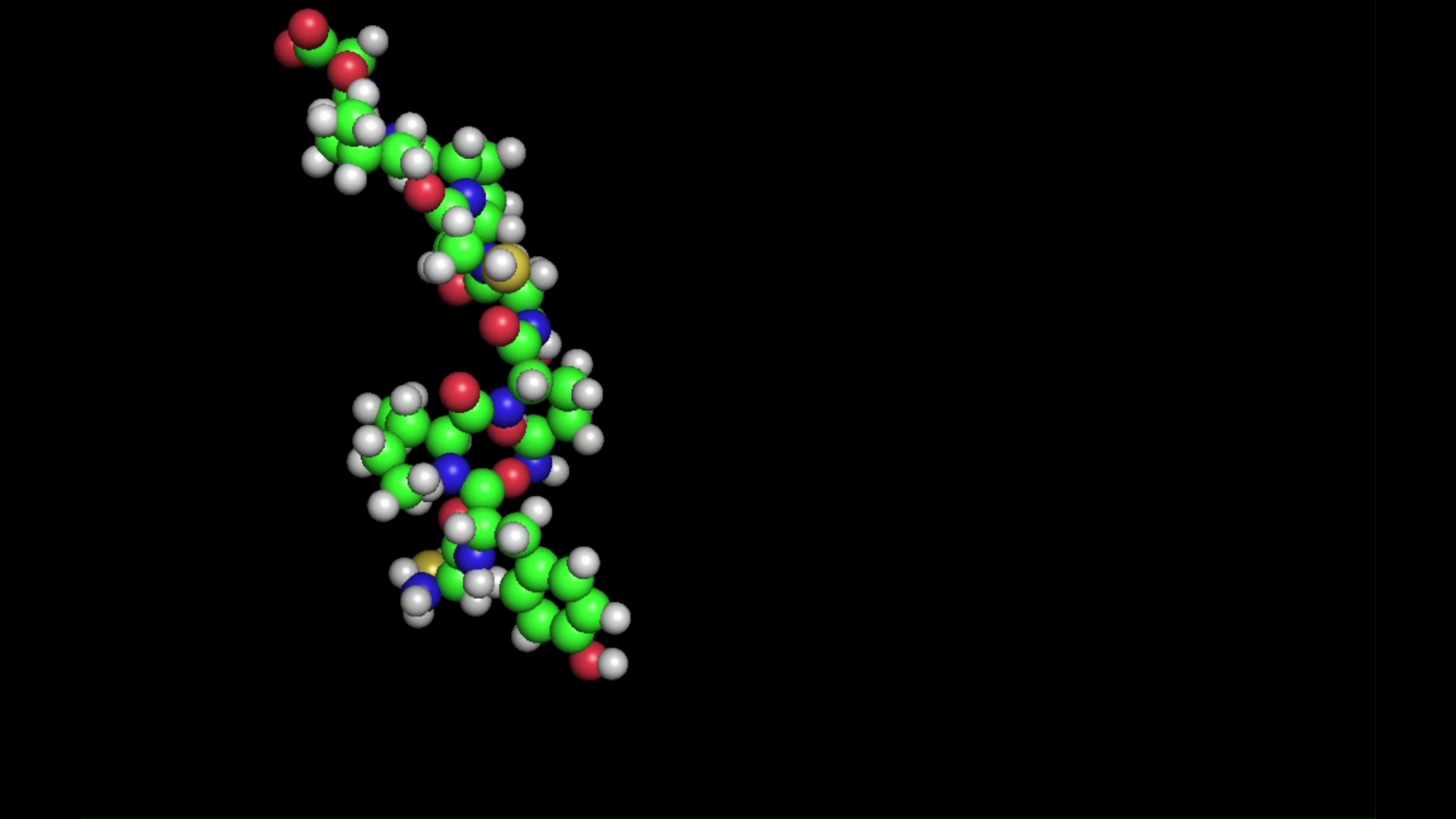}
\endminipage\hfill
\minipage{0.22\textwidth}%
  \includegraphics[width=\linewidth]{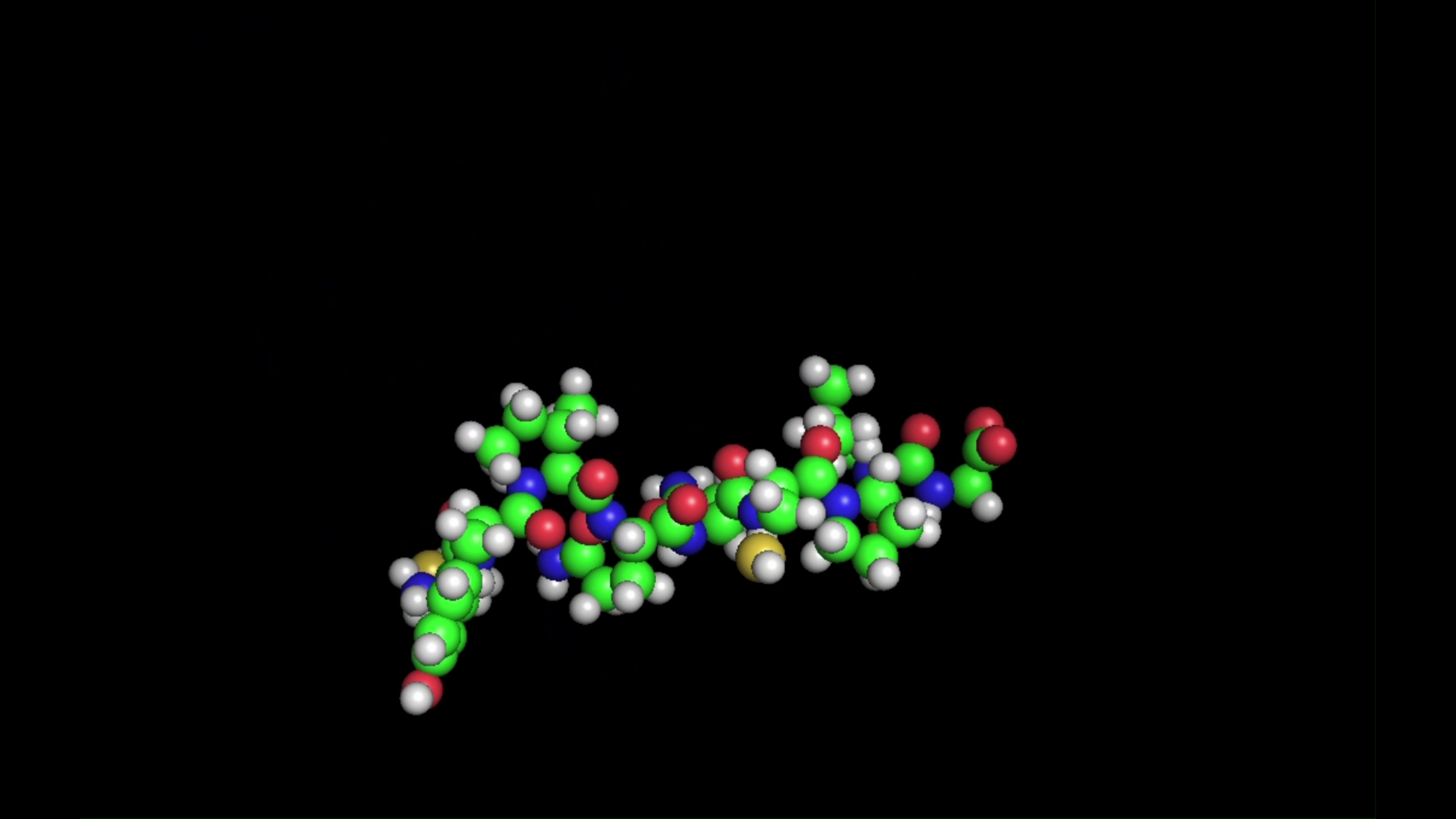}
\endminipage\hfill
\minipage{0.22\textwidth}%
  \includegraphics[width=\linewidth]{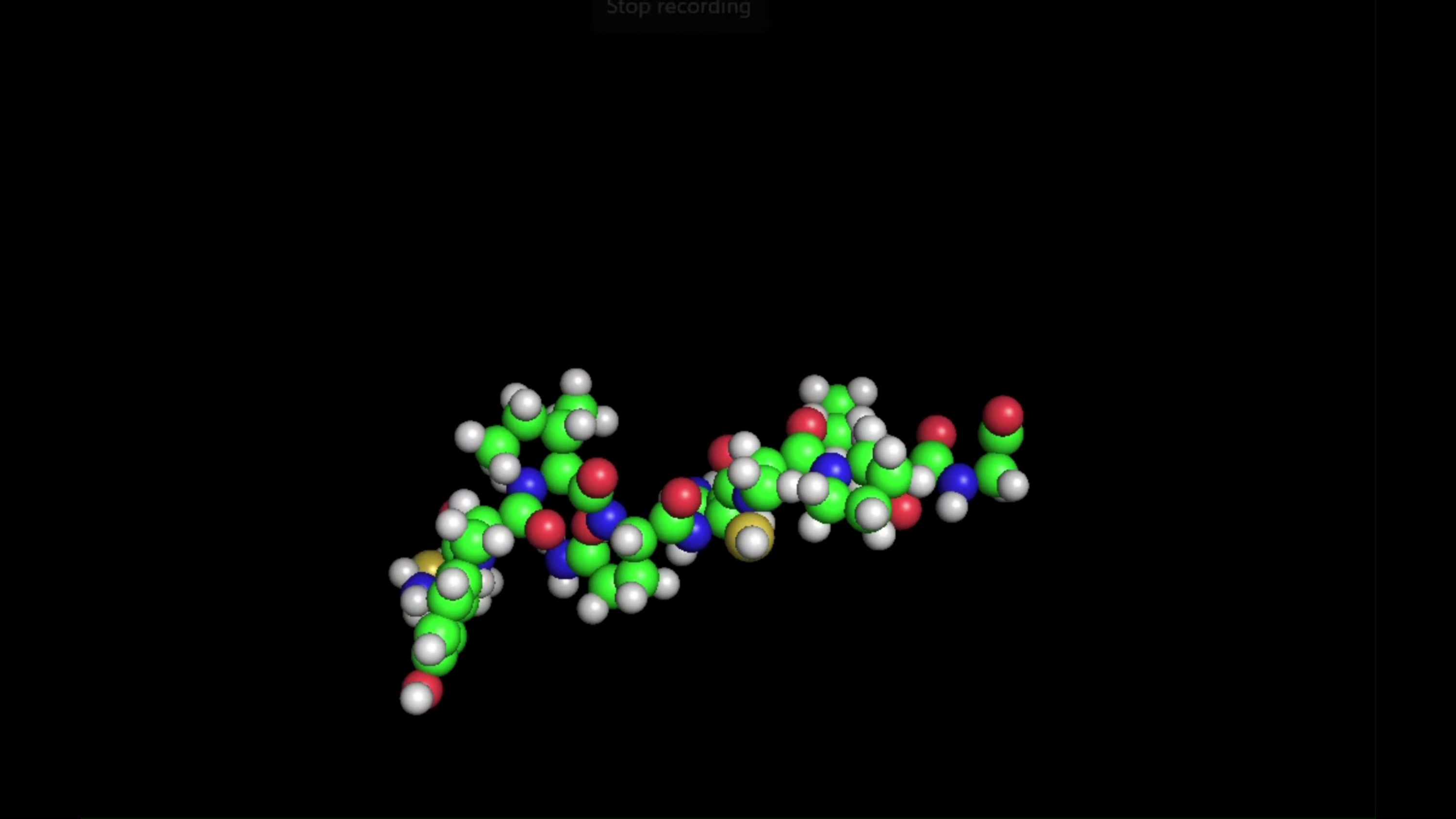}
\endminipage
\caption{Snapshots from a learned episode for Oxytocin molecule} \label{fig:visualize_oxytocin}
\end{figure}

\noindent For benchmarking purposes, we use the \textbf{Ratio} metric, which represents the proportion of episodes that successfully converge to stable, low-energy configurations, calculated from 100 episodes per molecule. An episode is considered successful if the system's energy consistently decreases throughout the trajectory, with the final energy falling within 10\% of the global minimum. Additionally, energy fluctuations in the final stabilization phase must remain within 10\% of the final energy value for at least the last 20\% of the episode. This metric quantifies the quality of the molecular trajectory by assessing both energy reduction and stability over time.\\

\noindent \cref{table:benchmarks} presents the benchmarking results, highlighting the performance of our method across six distinct molecules. Our method consistently outperforms both the Greedy and Random Hamiltonian approaches and shows competitive results against the NEMO-based method, especially for more complex molecules such as Oxytocin and Substance-P. While NEMO-based methods rely on Langevin dynamics for optimization, our reinforcement learning framework provides more flexibility and adaptability, especially in handling non-convex energy landscapes.\\

\noindent \cref{fig:energy_plot_over_time} shows the energy trajectory over 100 episodes for each molecule, plotted against time steps. These plots illustrate the overall stabilization process and the "butterfly" convergence pattern observed in our experiments. In these plots, we observe a distinctive butterfly-shaped convergence, where the molecules initially undergo significant fluctuations as the model explores the state space, followed by a gradual stabilization towards low-energy states. This convergence pattern results from the exploration-exploitation trade-off embedded in the reward function. During the stabilization phase, the model consistently narrows down the energy fluctuations, ensuring small variations around the optimal configuration.\\

\noindent The plots in \cref{fig:energy_plot_over_time} also illustrate the ability of our method to effectively escape local minima through strategic exploration, consistently achieving lower energy states than both the Greedy and Random Hamiltonian baselines. For example, in the case of Met-enkephalin (Figure 4(c)), one episode initially becomes trapped in a local minimum near an energy level of 1100 at approximately time step 120. However, the molecule later takes decisive high-magnitude actions, overcoming the local trap and eventually converging to a lower-energy state below 800. A similar scenario is observed for Bradykinin (Figure 4(a)), where the molecule becomes stuck at an energy level of around 4000 at time step 20, but successfully escapes by time step 50, reaching a significantly lower energy state below 2000.\\

\begin{figure}[!htb]
\minipage{0.22\textwidth}
  \includegraphics[width=\linewidth]{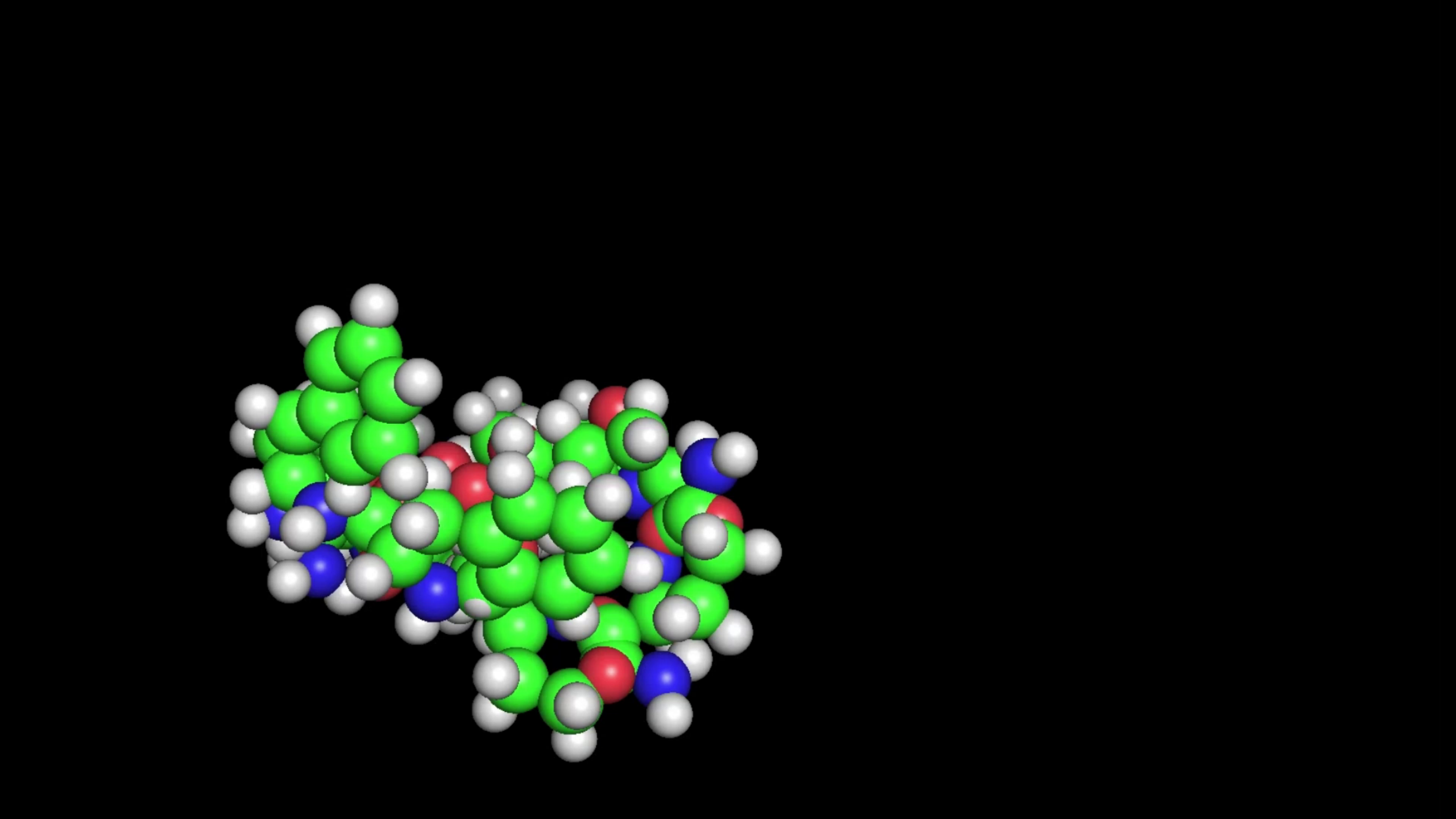}
\endminipage\hfill
\minipage{0.22\textwidth}
  \includegraphics[width=\linewidth]{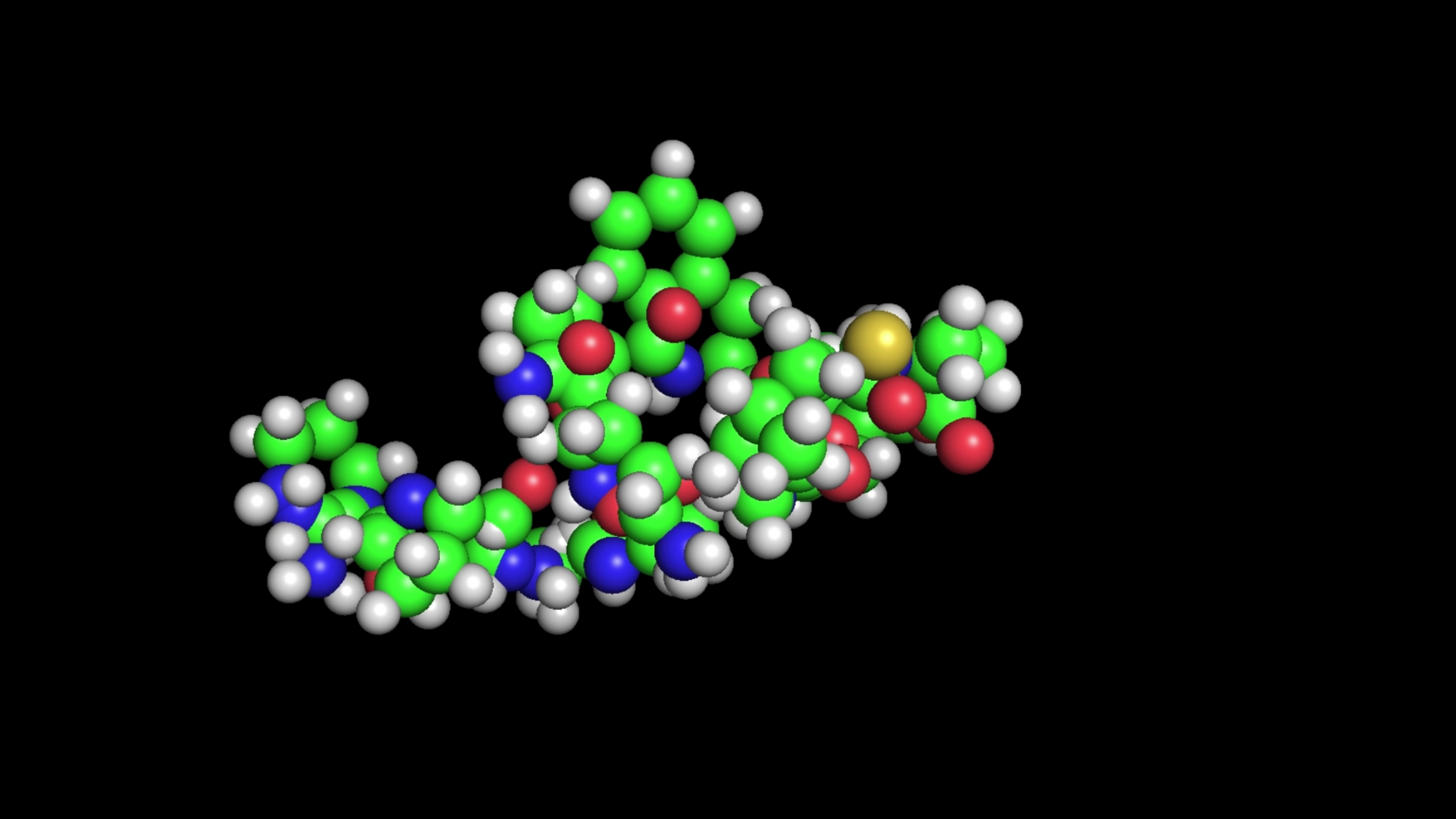}
\endminipage\hfill
\minipage{0.22\textwidth}%
  \includegraphics[width=\linewidth]{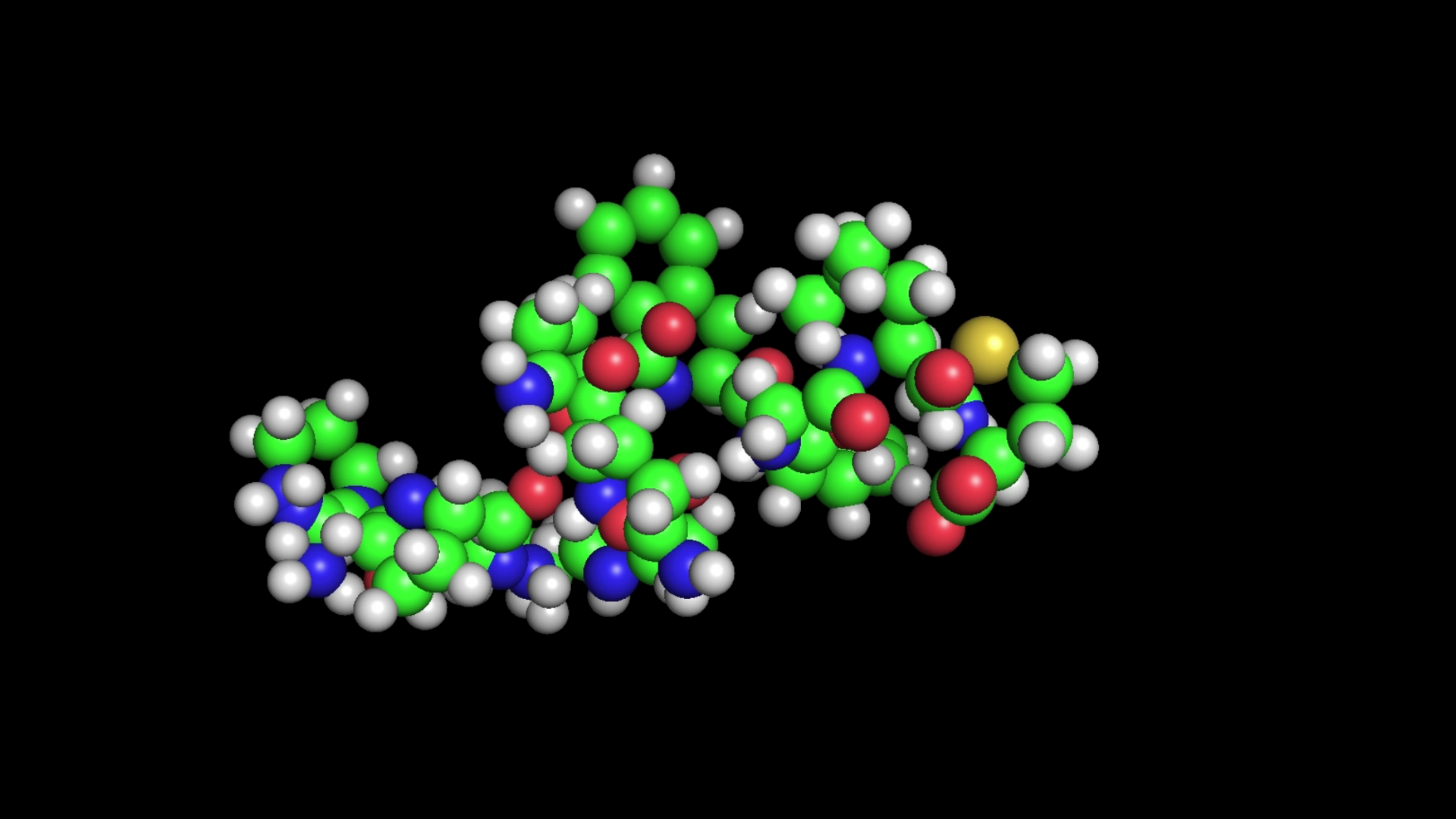}
\endminipage\hfill
\minipage{0.22\textwidth}%
  \includegraphics[width=\linewidth]{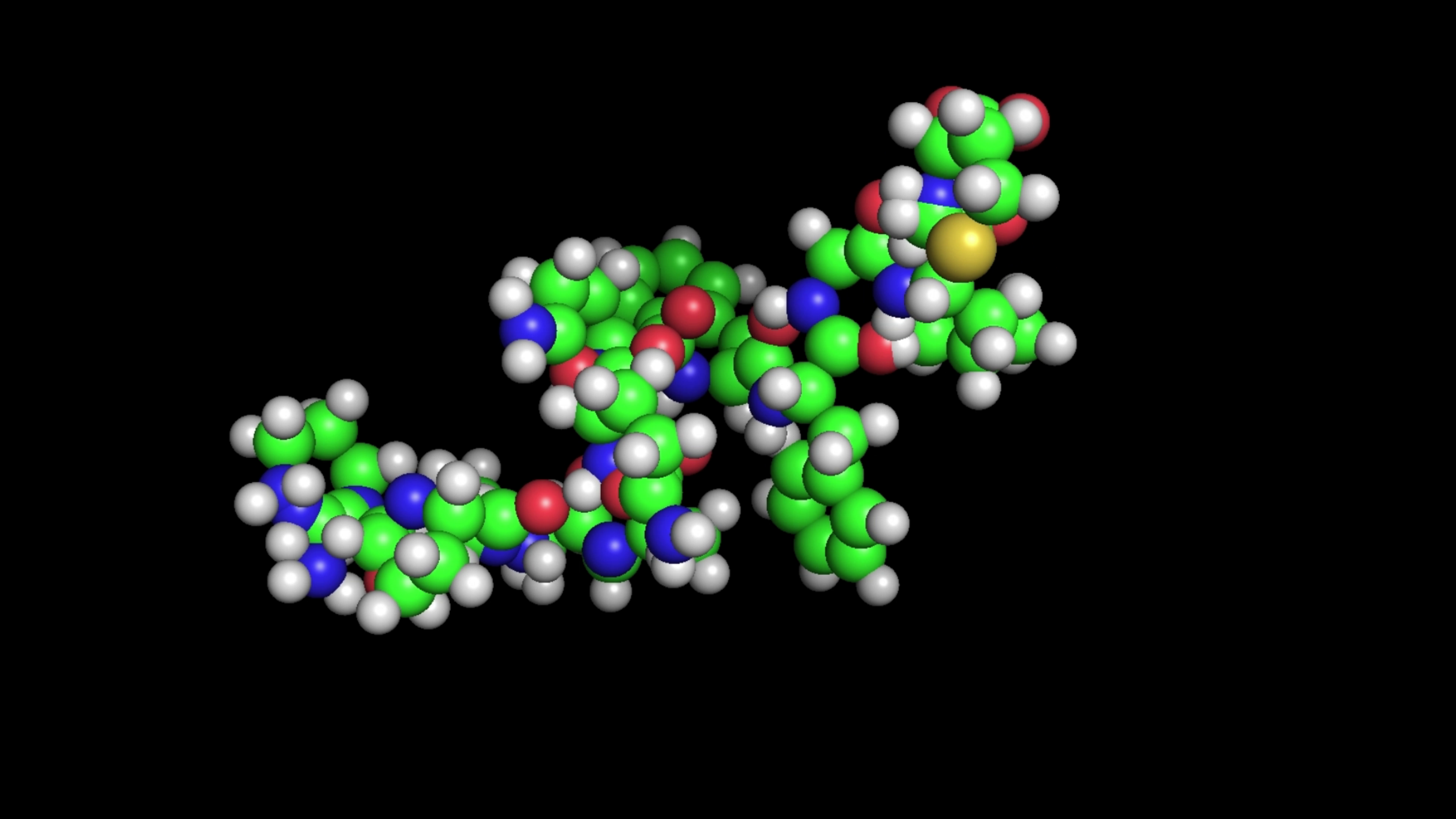}
\endminipage
\caption{Snapshots from a learned episode for Substance-P molecule} \label{fig:visualize_substance_p}
\end{figure}

\begin{figure}[!htb]
\minipage{0.22\textwidth}
  \includegraphics[width=\linewidth]{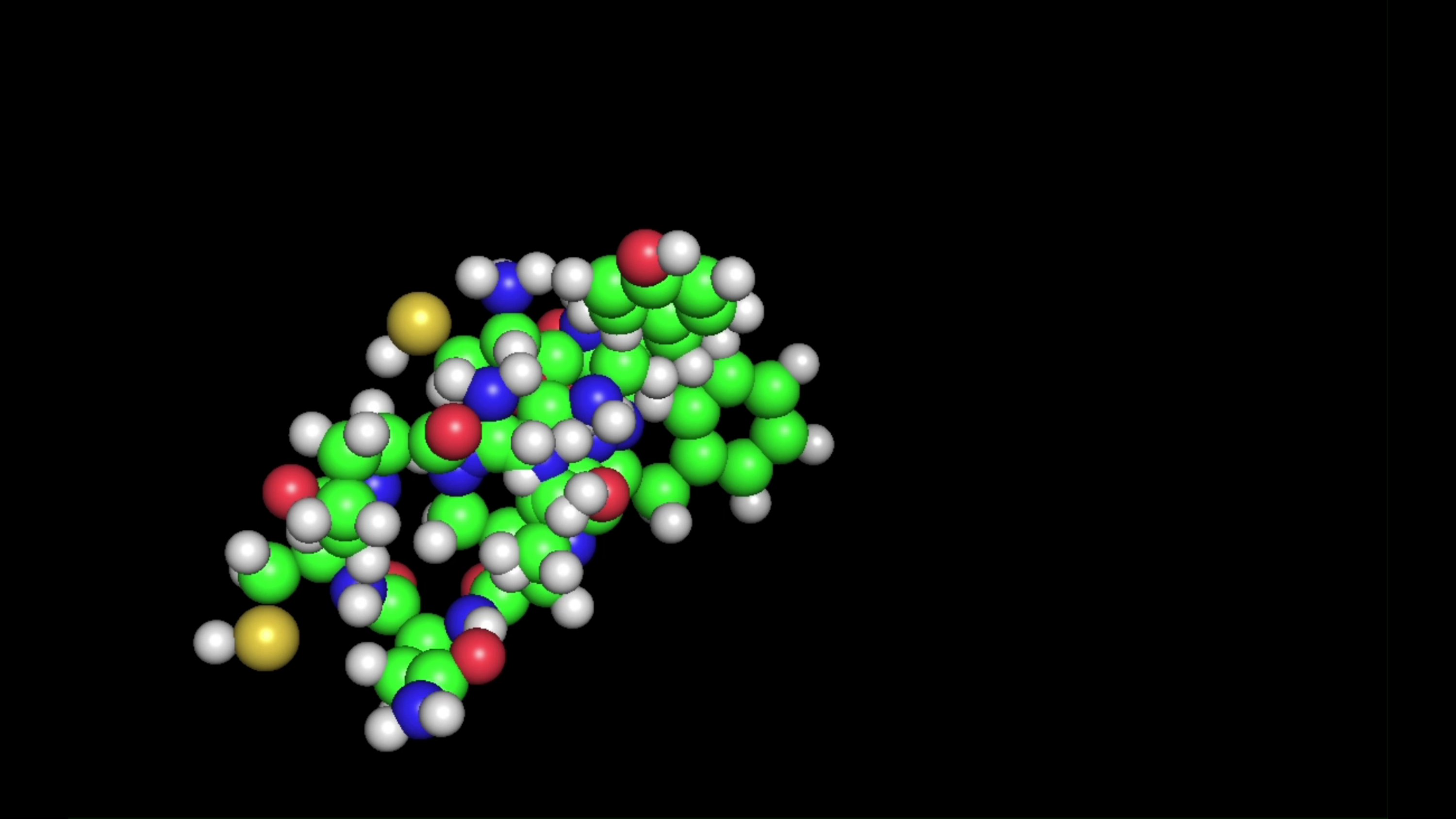}
\endminipage\hfill
\minipage{0.22\textwidth}
  \includegraphics[width=\linewidth]{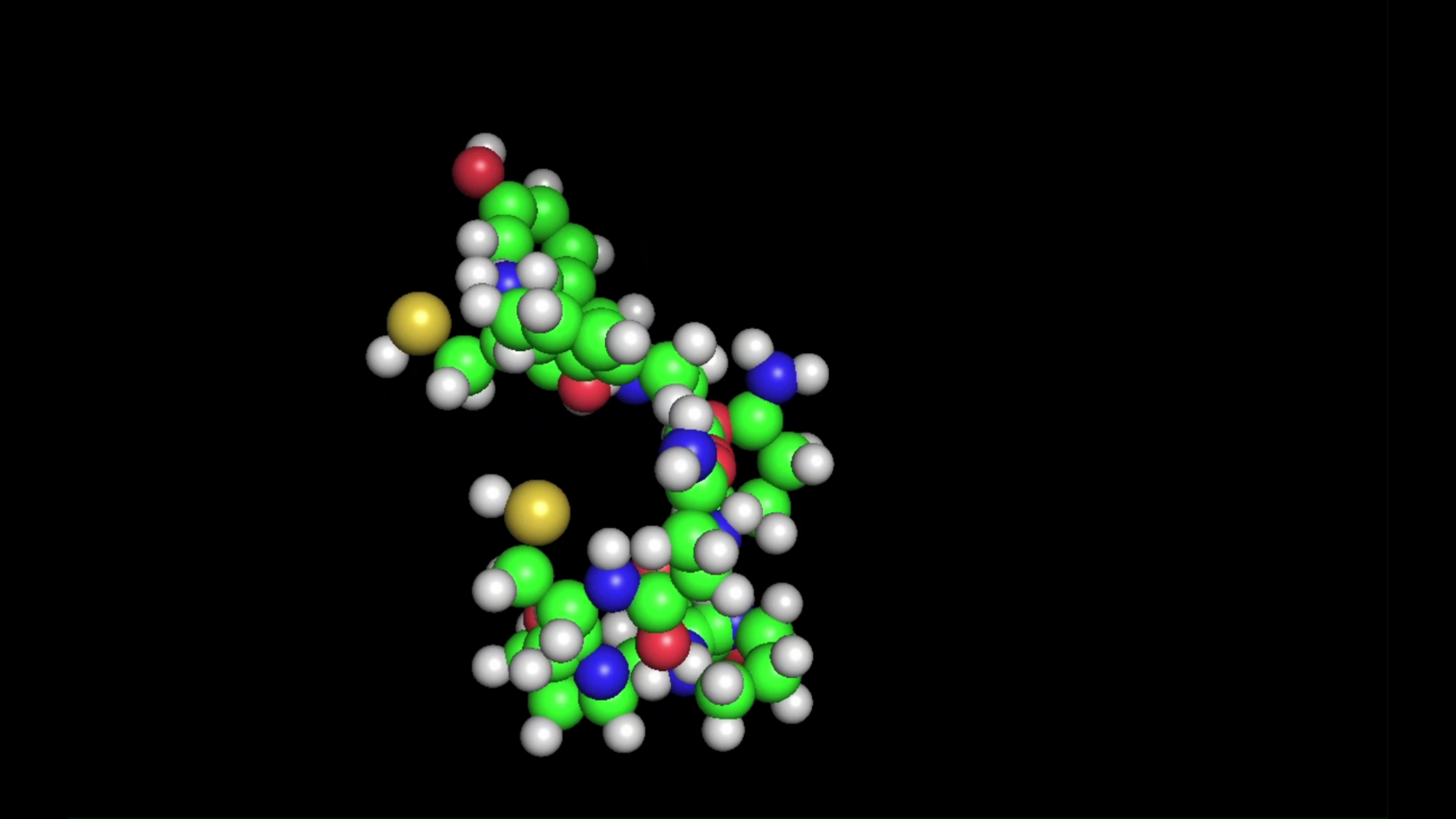}
\endminipage\hfill
\minipage{0.22\textwidth}%
  \includegraphics[width=\linewidth]{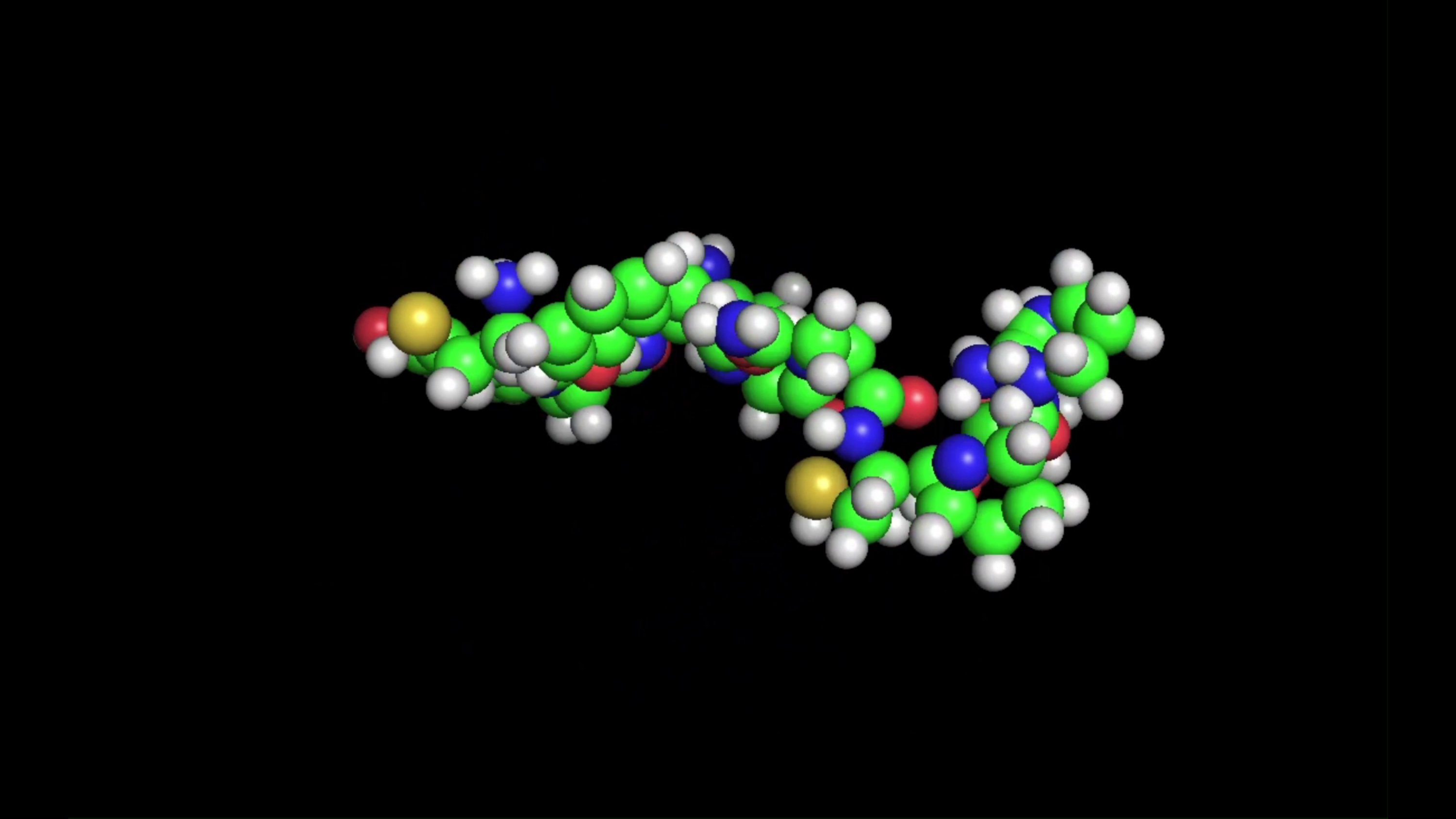}
\endminipage\hfill
\minipage{0.22\textwidth}%
  \includegraphics[width=\linewidth]{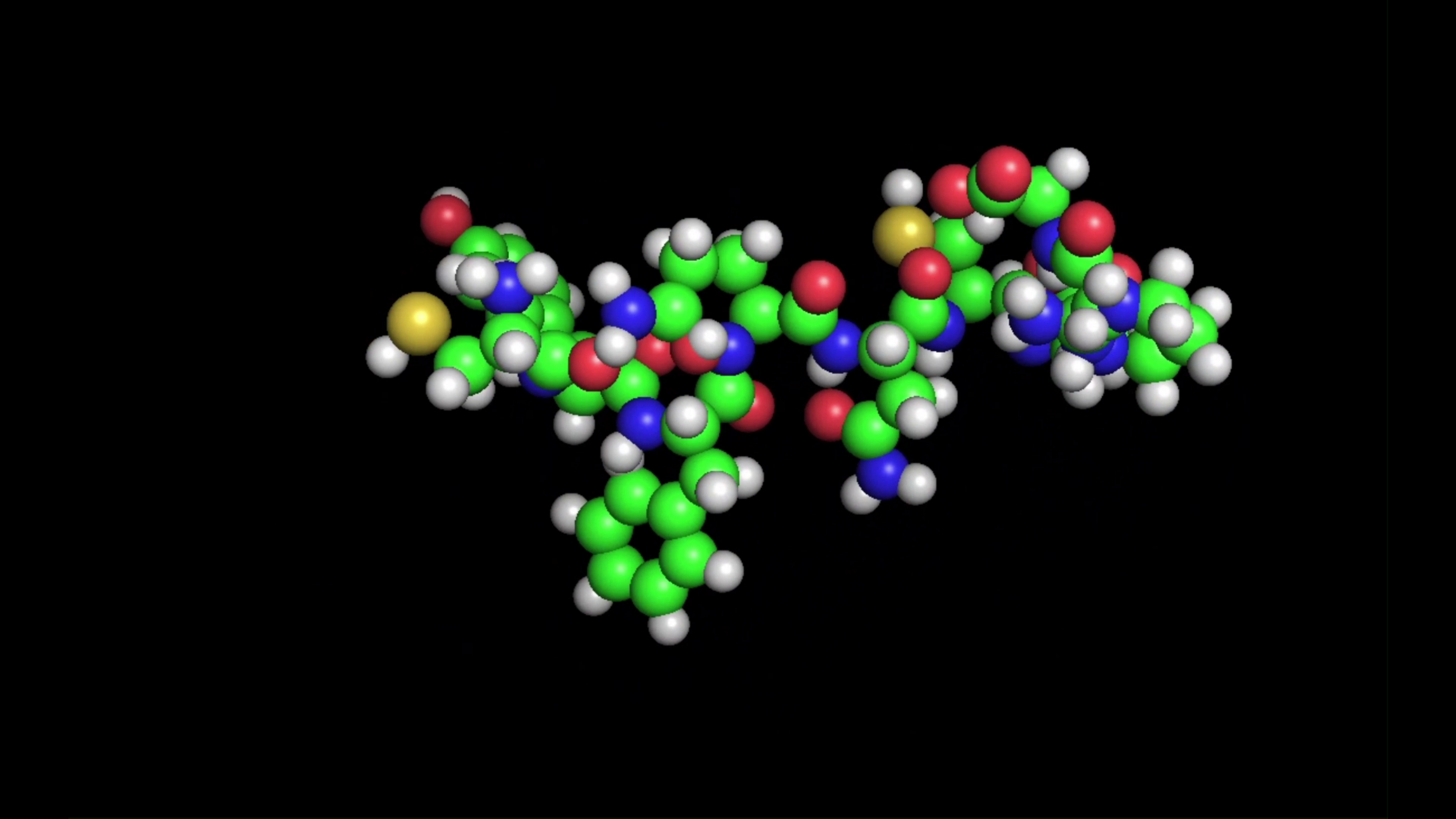}
\endminipage
\caption{Snapshots from a learned episode for Vasopressin molecule}\label{fig:visualize_vasopressin}
\end{figure}

\noindent To further understand the optimization process, we visualize the evolution of molecular structures at different time steps. \cref{fig:visualize_bradykinin}, \cref{fig:visualize_cln025}, \cref{fig:visualize_met_enkephalin}, \cref{fig:visualize_oxytocin}, \cref{fig:visualize_substance_p}, and \cref{fig:visualize_vasopressin} show snapshots of each molecule during various phases of the learning process. These visualizations provide insights into how the molecular configurations shift dynamically as the reinforcement learning model optimizes the trajectory. In these figures, we observe how the molecular configurations evolve from initial high-energy states to more stabilized configurations over time. The dynamic trajectory is captured throughout the learning process, showcasing continuous structural adjustments rather than abrupt changes. This gradual shift is a key advantage of our approach. It maintains control over the molecular system’s evolution, minimizing large energy jumps and ensuring smooth transitions towards optimal configurations.\\

\noindent For example, in \cref{fig:visualize_met_enkephalin} (Met-enkephalin), the molecule initially adopts a loosely packed configuration with large energy gaps. As the model progresses through learning episodes, the molecule gradually packs into a more stable conformation, minimizing both energy and structural variance. Similarly, \cref{fig:visualize_vasopressin} (Vasopressin) illustrates how the molecule transitions from a highly fluctuating state to a tightly bound, low-energy configuration, stabilizing in the later stages of the optimization.

\section{Conclusions}\label{sec:conclusion}
In this paper, we introduced a reinforcement learning framework for optimizing molecular dynamics, utilizing a continuous dynamical approach based on the stochastic Pontryagin Maximum Principle (PMP). Our method focuses on optimizing the entire molecular trajectory, enabling efficient exploration of complex energy landscapes, escaping local minima, and converging on stable low-energy states. By integrating the Soft Actor-Critic (SAC) algorithm, our framework operates in an unsupervised manner, requiring no labeled data, and provides precise control over molecular behavior during both exploration and stabilization phases. Through comprehensive experiments on various molecules, we demonstrated that our approach achieves superior performance compared to other unsupervised physics-based methods. This trajectory-based optimization method is well-suited for practical applications, such as optimizing drug binding affinities, improving molecular designs for material science, or assisting in computational simulations for personalized medicine. Future work will focus on scaling our approach to address larger molecular systems with intricate interactions and enhancing computational efficiency for multi-scale simulations.\\

\noindent The implementation is available at: \textcolor{red}{\href{https://github.com/CVC-Lab/SAC-for-H-Bond-Learning}{https://github.com/CVC-Lab/SAC-for-H-Bond-Learning}}.

\begin{credits}
\subsubsection{\ackname} This research was supported in part by a grant from NIH- R01GM117594, by the Peter O’Donnell Foundation and in part from a grant from the Army Research Office accomplished under Cooperative Agreement Number W911NF-19-2-0333.\\

\noindent This preprint has no post-submission improvements or corrections. The Version of Record of this contribution is published in the Neural Information Processing, ICONIP 2024 Proceedings.
\end{credits}

%
%
%
\bibliographystyle{splncs04}
\bibliography{citation}

\appendix
\section{Proof of \cref{lem:dynamic_form}}
To prove the lemma, we begin by revisiting the continuous-time stochastic control formulation with the following objective:
\begin{equation}\label{control_prob}
J(a(.)) = \E \left[\int_0^T l(s_t, a_t)dt\right]
\end{equation}
subject to the randomized dynamics:
\begin{align}
ds_t &= g(s_t, a_t)dt + \sigma(s_t, a_t)dW_t\\
s(0) &= s_0
\end{align}
Here, $s_t$ is the state variable at time $t$ with values in $\mathbb{R}^n$, $a(.)$ is the control/action variable with values in a subset of $\mathbb{R}^n$, and $W_t$ is the Brownian motion in $\mathbb{R}^n$. Additionally, we consider the following Hamiltonian function:
\begin{equation}\label{Hamiltonian}
H(s, a, p, K) = l(s, a) + (p, g(s, a)) + \sum_{j=1}^d (K_j, \sigma^j(s, a)) 
\end{equation}

\noindent Let the first-order adjoint process $(p(.), K(.))$ be the unique process that satisfies the stochastic differential equation:
\begin{align}
-dp_t &= \bigg(g_s(s_t, a_t)^T p_t + \sum_{j=1}^d (\sigma_s^{j}(s_t, a_t))^T K_j(t) + l_s(s_t, a_t)\bigg)dt - K_tdW_t\\
p_t &= h_s(s(T))
\end{align}

\noindent Using the result in \cite{Peng1990-gc}, under continuously differentiable constraints, we obtain the maximum principle $H_a(s_t, a_t, p_t, K_t) = 0$ for the optimal path $(s_t, a_t)_{t = 0}^T$. Now, applying this maximum principle to the special case of our control problem \cref{control_prob} where the covariance $\sigma$ is constant, with drift $g(s_t, a_t) = a_t$, and the running loss $l(s_t, a_t) = -r(s_t, a_t) = \Gamma(t)(V(s_t) \mp a_t^2/2)$, we find that $0 = H_a(s_t, a_t, p_t, K_t) = p_t + l_a(s_t, a_t)$. Hence $p_t = -l_a(s_t, a_t) = \pm \Gamma(t) a_t$. For this special case, we also derive the following dynamics equation for the adjoint variable $p_t$:
\begin{equation}
-dp_t = 0 + 0 + l_s(s_t, a_t) dt - K_t dW_t
\end{equation}

\noindent Now we have:
\begin{equation}
\pm(\Gamma(t) da_t + \Gamma'(t)a_tdt) = \pm \frac{d}{dt}(\Gamma(t) a_t)dt = dp_t = -l_s(s_t, a_t) dt + K_t dW_t
\end{equation}

\noindent This simplifies to:
\begin{equation}
\Gamma(t) da_t + \Gamma'(t) a dt = \mp \Gamma(t)\nabla V(s_t) dt + \pm K_t dW_t
\end{equation}

\noindent As a result, the final dynamic equation for the optimal $a_t$ is:
\begin{equation}
da_t = -\frac{\Gamma'(t)}{\Gamma(t)} a_t dt + \mp \nabla V(s_t) dt + \pm \frac{K_t}{\Gamma(t)} dW_t
\end{equation}
 
\noindent Due to the uniqueness of the solution to the adjoint stochastic differential equations and the uniqueness of the SDE in the form of the Langevin dynamics, $\pm \frac{K_t}{\Gamma(t)} = K$, and we obtain the final dynamics as desired:
\begin{align}
ds_t &= a_t dt + \sigma dW_t \\
da_t &= -\frac{\Gamma'(t)}{\Gamma(t)} a_t dt \mp \nabla V(s_t) dt + K dW_t
\end{align}

\end{document}